\documentclass[smallcondensed,final]{svjour3}                     
\smartqed  

\usepackage{amsmath}
\usepackage{graphicx}
\usepackage{amssymb}
\usepackage{xcolor}
\usepackage{lmodern} 

\newcommand{\mb}[1]{{\mathbf{#1}}}
\newcommand{\ds}{\displaystyle}

\newcommand{\dd}[2]{{\partial_{#1} #2}}

\newcommand{\be}{\begin{equation}}
\newcommand{\eb}{\end{equation}}
    
\title{Fourier Continuation Discontinuous Galerkin Methods for Linear Hyperbolic Problems\thanks{This work  
was supported by NSF Grant DMS-1913076. Any opinions, findings, and conclusions or recommendations expressed in this material are those of the authors and do not necessarily reflect the views of the National Science Foundation}}

\titlerunning{FC-DG Methods for Linear Hyperbolic Problems}        
\author{Kiera van der Sande  \and 
Daniel Appel\"{o} \and 
Nathan Albin
}
\institute{Kiera van der Sande\at
  Department of Applied Mathematics, University of Colorado Boulder, Boulder, CO, USA. \\
              \email{kiera.vandersande@colorado.edu}
           \and
Daniel Appel\"{o} \at Department of Computational Mathematics, Science \& Engineering and Department of Mathematics, Michigan State University, East Lansing, USA.  \\
\email{appeloda@msu.edu}
           \and
Nathan Albin \at Department of Mathematics, Kansas State University, Manhattan, KS, USA. \\
\email{albin@k-state.edu}
}

\date{Received: date / Accepted: date}

\newcommand{\myfig}{Figure }

\begin{document}
	\maketitle
	
\begin{abstract}
Fourier continuation is an approach used to create periodic extensions of non-periodic functions in order to obtain highly-accurate Fourier expansions. These methods have been used in PDE-solvers and have demonstrated high-order convergence and spectrally accurate dispersion relations in numerical experiments. Discontinuous Galerkin (DG) methods are increasingly used for solving PDEs and, as all Galerkin formulations, come with a strong framework for proving stability and convergence. Here we propose the use of Fourier continuation in forming a new basis for the DG framework. 
\end{abstract}

\section{Introduction}
When approximating solutions to partial differential equations the choice of functions to use in the approximation impacts the accuracy, efficiency and stability of the resulting numerical method. For time dependent wave propagation problems on bounded domains most methods use a polynomial approximation. This can be done through local polynomials that interpolate discrete function values at grid-points as is done in finite difference methods \cite{GusKreOli95}. Numerical derivatives are then obtained by analytic differentiation of the interpolant. Another approach is Galerkin's method, which starts from the variational formulation of the equations and seeks a polynomial approximation such that the residual of the approximated PDE is orthogonal to all polynomials in the approximation space, \cite{ern2013theory}. In particular for wave propagation problems the discontinuous Galerkin (DG) method \cite{cockburn1989tvb,cockburn2001runge,Hesthaven:2008fk} has emerged as an accurate and robust approach. However, the high degree polynomial approximation on each element that is used in DG and spectral elements results in numerical stiffness and reduces the allowable timesteps significantly below the limit dictated by physical finite speed of propagation considerations. This can limit the efficiency, particularly for linear hyperbolic systems of equations.      

It is widely known that periodic functions are well approximated by Fourier series or trigonometric interpolation, and that these approximations on a uniform grid can be computed and manipulated efficiently using discrete fast Fourier transforms (FFTs). Non-periodic functions that are sampled on a grid may still be approximated by trigonometric interpolation, but the approximation becomes oscillatory and inaccurate near boundaries due to Gibbs' phenomenon. There has been much interest in overcoming this problem, including the approach known as Fourier continuation (FC) wherein a periodic extension allows non-periodic functions to be represented as a trigonometric series. Several FC methods, also known as Fourier extension, have been developed and have shown superalgebraic and even exponentially accurate approximation properties \cite{HuybrechsFE2010,LyonFastFC2011,boydFEcomparison2002}. In general, they seek a Fourier series representation which is close in the least-squares sense to the original function on a bounded interval.
	
Fourier continuation methods, particularly the FC-Gram approach, have been used in several PDE-solvers where they have demonstrated high order convergence rates combined with very small dispersive errors \cite{BrunoFCADBasic2010,LyonFCADII2010,BrunoFCADVarCoeff2012,AlbinFCNavier2011,AlbinAcoustic2012}. In addition, the approximations to derivatives obtained by the FC approach cause considerably less numerical stiffness than those of DG. However, although successful for many applications, to our knowledge, FC-based numerical PDE solvers does not come with a provable guarantee of stability. The method we propose here is, in its current incarnation, not as fast as previous FC-PDE solvers but it does come with the usual stability guarantee intrinsic to Galerkin formulations. And, being an element based discretization, it can handle geometry by the use of unstructured meshes.      
	
Given the solid theoretical foundation and robustness of the discontinuous Galerkin method and the small dispersive errors and large timesteps of FC-based PDE solvers it is natural to combine the two. In this paper we propose a new basis constructed using Fourier continuation to create functions that are periodic on an extended domain. We then use this basis for constructing DG discretizations for linear hyperbolic equations such as transport equations and  Maxwell's equations. 

As we show though numerical experiments the resulting FC discontinuous Galerkin methods have small spectral radius, allowing large timesteps, and their dispersive properties results in methods that can propagate waves over long distances with minimal dispersive errors. A drawback of the FC-DG method is that in general the FC basis will not be orthogonal, which leads to dense (but well conditioned) mass matrices and stiffness matrices. Here we are mainly concerned with the approximation properties of the method and delay efficient implementations to the future. We note that rapid inversion of the mass matrix and application of the stiffness matrix will likely require us to adopt matrix free approaches such as those in \cite{KroKor19}. We also note that Bruno and Prieto,  \cite{bruno2014spatially}, have demonstrated that variable coefficient boundary value problems discretized by FC methods can be solved very efficiently by finite difference preconditioned GMRES.

The rest of this paper will be organized as follows. In Section \ref{sec:FCbasis}, we explain the Fourier continuation method and how it is used to generate a basis. In Section \ref{sec:DGimplementation}, the DG formulation is reviewed and we explain the methods we will use for solving PDEs with our proposed basis, including considerations for computing integrals and time-stepping. Section \ref{sec:numericalexperiments} contains numerical experiments and results for test problems in 1-D and 2-D, and Section \ref{sec:EMwaves} considers applications to electromagnetic waves for problems such as in optical media. Conclusions and future directions are outlined in Section \ref{sec:conclusion}.
	
\section{Fourier Continuation as a Basis \label{sec:FCbasis}}

\begin{figure}
\includegraphics[width=\textwidth]{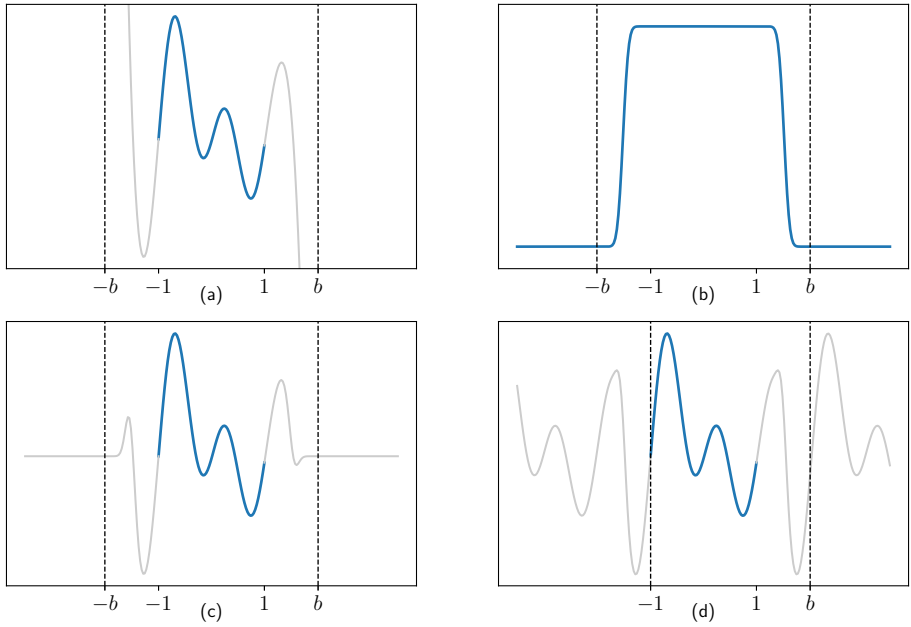}
\caption{Steps of the Fourier continuation construction of a discrete periodic extension.}
\label{fig:fc-steps}
\end{figure}

Given the optimal approximation properties of Fourier expansions and the speed at which they can be manipulated, we use a version of the the FC(Gram) Fourier continuation method introduced in~\cite{BrunoFCADBasic2010,LyonFastFC2011} to generate a basis for DG methods.  Specifically, we use the modified construction, developed in~\cite{AlbinDPE2014}, for generating a discrete periodic extension of a smooth function sampled on a uniform grid.

Conceptually, the periodic extension operator can be visualized through the series of steps shown in \myfig~\ref{fig:fc-steps}.  The construction begins with the values $f_l=f(z_l)$ of a smooth function on the interval $[-1,1]$ sampled on the uniform grid $z_l=-1 + 2l/(N-1)$ for $l=0,1,\ldots,N-1$.  These values are indicated by the dark curve in \myfig~\ref{fig:fc-steps}a.  The values are then extended to a larger interval, $[-b,b]$ using polynomial extrapolation.  More precisely, the construction depends on two positive integer parameters, $p$ and $M$.  First, the function is extended to the grid points $z_l$ with $l=-M,-M+1,\ldots,-1$ using the $(p-1)$-degree polynomial interpolant of the left-most $p$ samples.  Similarly, the function is extended rightward to the grid points $z_l$ with $l=N,N+1,\ldots,N+M-1$ using the polynomial interpolant of the right-most $p$ samples.  This extends the samples to the interval $[-b,b]$ with $b=1+2M/(N-1)$.  This extension is indicated by the light curve in \myfig~\ref{fig:fc-steps}a.

Next, the extrapolated extension is multiplied by a smooth window function, shown in \myfig~\ref{fig:fc-steps}b.  This window function was constructed in~\cite{AlbinDPE2014} to have rapidly decaying Fourier coefficients and to be well-resolved on the discrete grid.  Moreover, with an error on the order of machine epsilon, the window function equals $1$ on $[-1,1]$ and $0$ outside $[-b,b]$.  The result of multiplying the function in \myfig~\ref{fig:fc-steps}a with the window function in \myfig~\ref{fig:fc-steps}b is shown in \myfig~\ref{fig:fc-steps}c.  In this way, we have extended the original sample values, $f_l$ to samples $\tilde{f}_l$ on a uniform grid on the whole real line.  By construction, these samples agree (up to a small error on the order of machine epsilon) with the original samples on the points $z_l$ for $l=0,1,\ldots,N$.  Moreover, the samples $\tilde{f}_l=0$ for $z_l$ with $l<-M$ or $l\ge N+M$.

To complete the periodic extension, we define the values
\begin{equation*}
f^c_l = \sum_{r=-\infty}^\infty \tilde{f}_{l+r(N+M)}\quad\text{for }l\in\mathbb{Z}.
\end{equation*}
The result is a discrete periodic function satisfying $f^c_l=f_l$ for $l=0,1,\ldots,N-1$ (to machine precision).  These values can be viewed as samples of a smooth, periodic function on the interval $[-1,b]$, as shown in \myfig~\ref{fig:fc-steps}d.

As described in~\cite{BrunoFCADBasic2010,LyonFCADII2010,AlbinDPE2014}, this procedure can be accelerated by pre-computing a linear extension operator mapping the values $f_l$ (actually only, the first $p$ and last $p$ values) to the extension values $f^c_l$ for $l=N,N+1,\ldots,N+M-1$.  Using the FFT, we can then find coefficients $a_k$ of a trigonometric polynomial
\begin{equation*}
f^c(z) = \sum_{k = -W}^W a_k \exp\left(\frac{2\pi i k z}{1+b}\right),
\end{equation*}
with the property that where $f^c(z_l) = f_l$ for $l = 0,1,...,N-1$.  The $W$ in the formula can be taken to be $W=\lfloor (N+M-1)/2\rfloor$.  Provided the samples $f_l$ came from a smooth, sufficiently resolved function on $[-1,1]$, the function $f^c$ will approximate $f$ on the same interval with high accuracy.

In order to produce a basis for the DG method, we apply the discrete periodic extension operator to the canonical basis $\lbrace \textbf{e}_i\rbrace_{i=1}^N$, which allows the basis to be represented in terms of its Fourier coefficients. Differentiation and spectral interpolation of the basis can then be done efficiently using the FFT.  The basis depends on a number of parameters including the number of discretization points $N$, the polynomial approximation order $p$, and the extension length $M$.  In particular, the parameter $p$ directly affects the order of accuracy of the method. Unless otherwise specified, we will use 9th degree interpolating polynomials in the Fourier continuation (i.e. $p = 10$ points) and $M = 25$ points in the extended domain, as in \cite{AlbinDPE2014}.

\begin{figure}[]
\begin{center}
\includegraphics[width=0.45\textwidth]{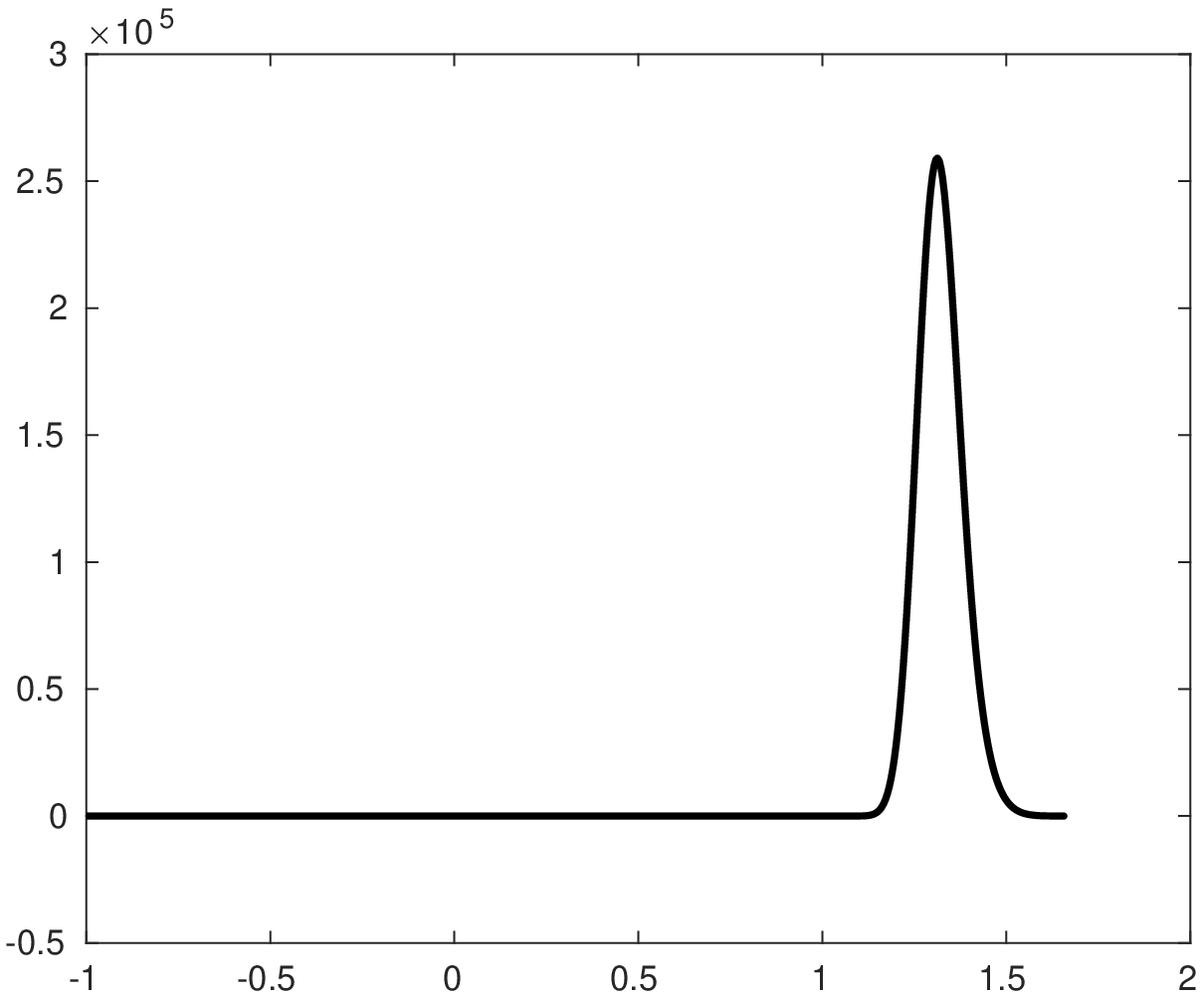}
\includegraphics[width=0.45\textwidth]{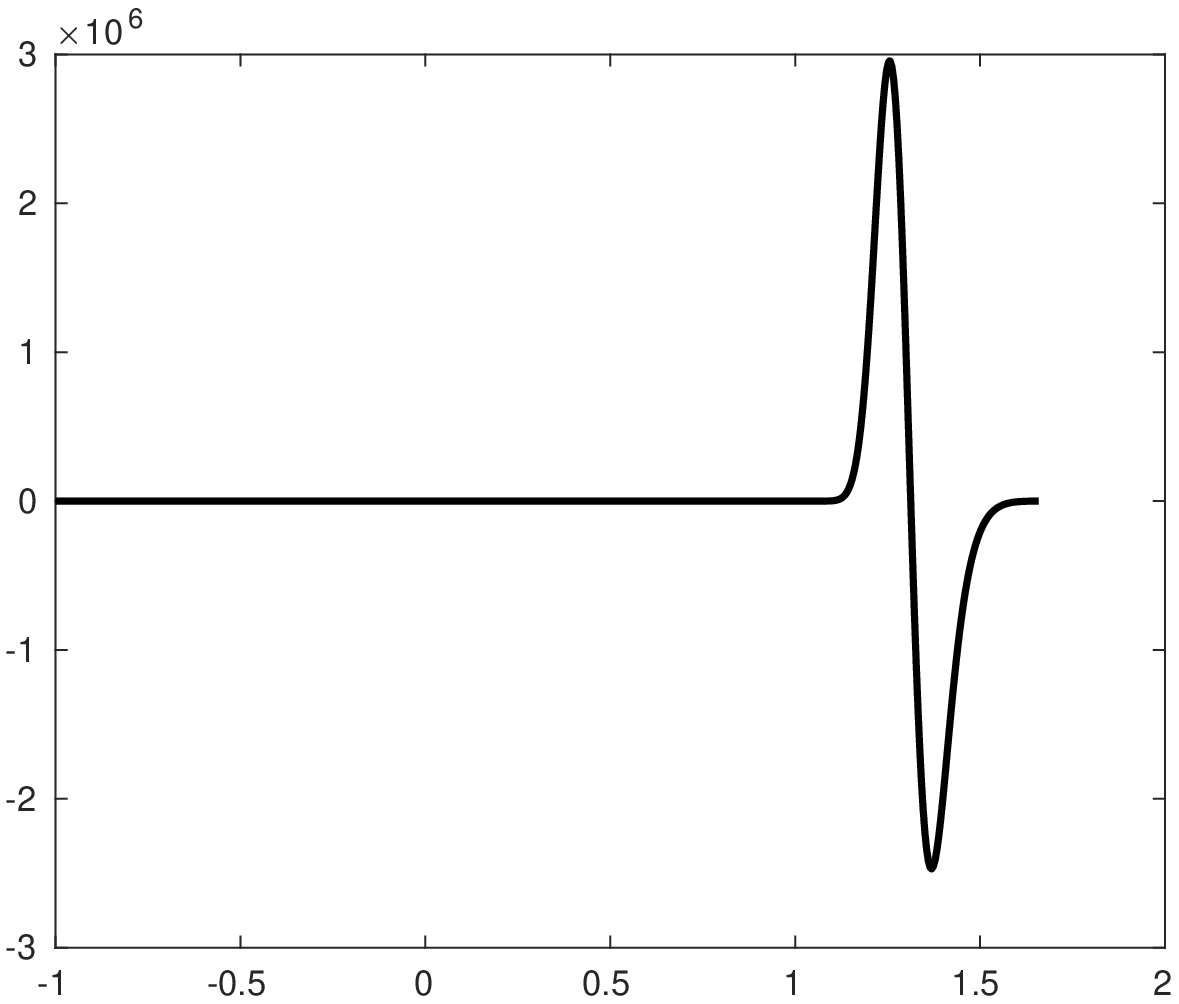}
\includegraphics[width=0.45\textwidth]{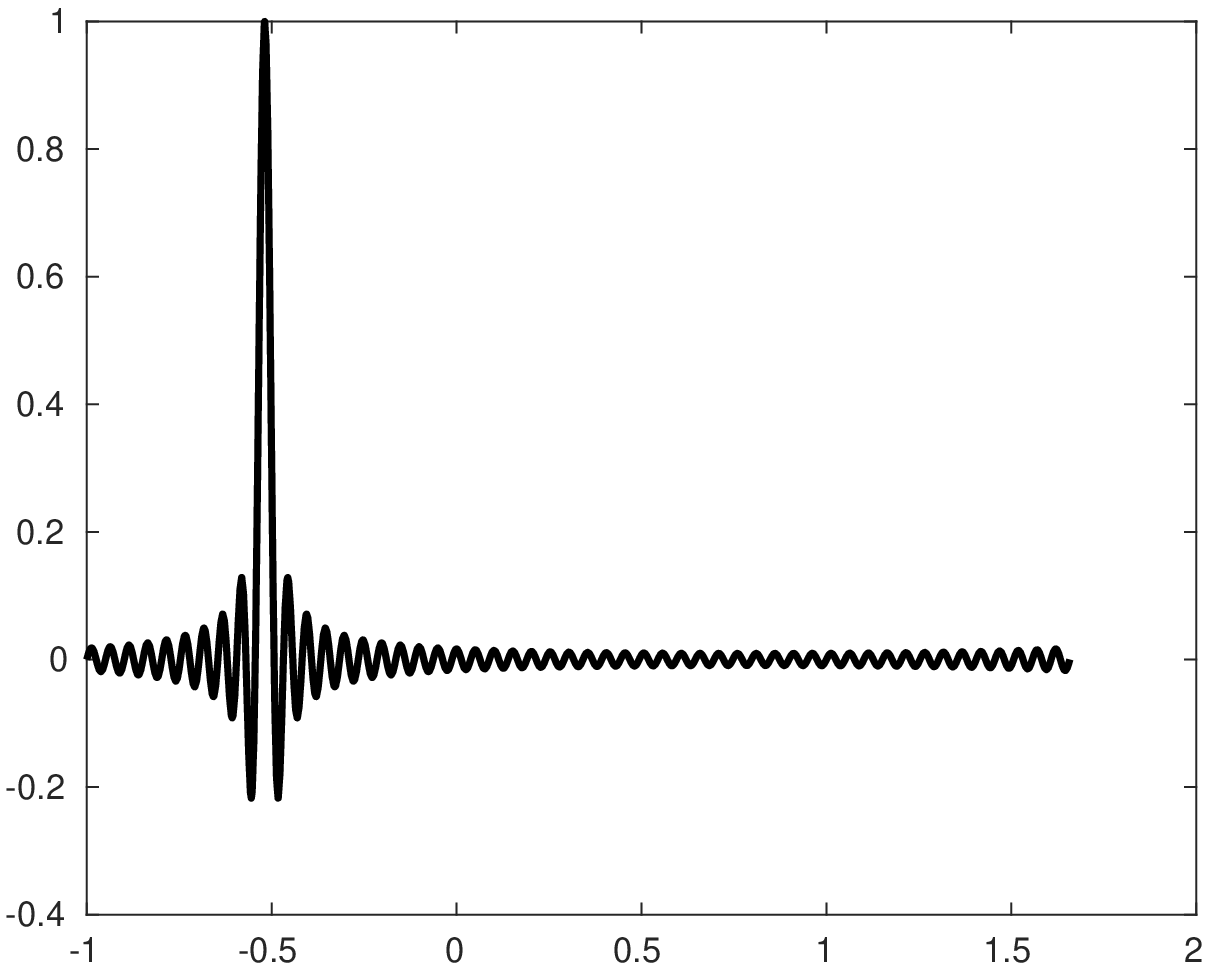}
\includegraphics[width=0.45\textwidth]{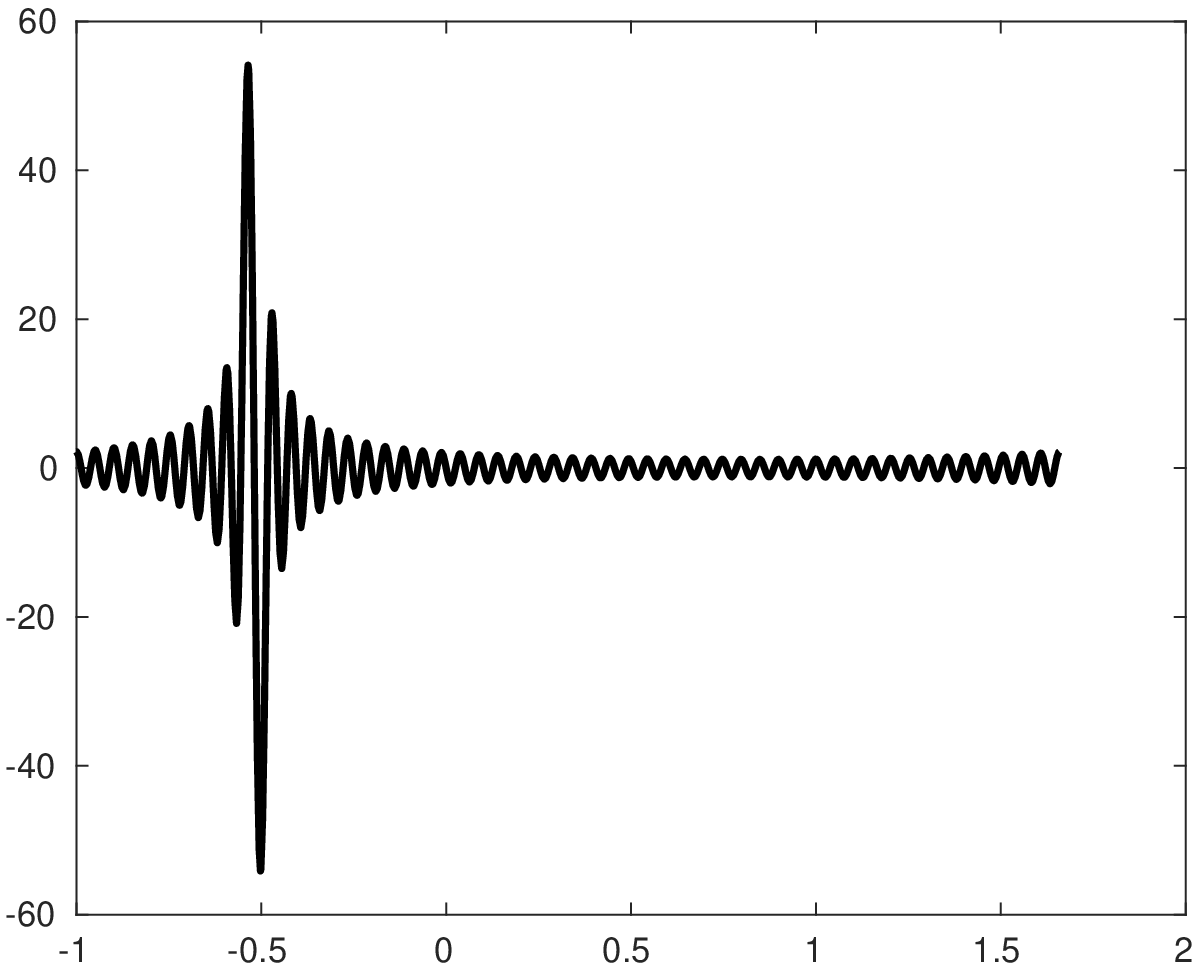}
\caption{Two basis functions, $\phi_1$ (top left) and $\phi_{20}$ (bottom left), resulting from the Fourier continuation of ${\bf e}_1$ and ${\bf e}_{20}$, and their derivatives (top / bottom right) on the extended periodic domain with $N = 80$ points on $[-1,1]$ and $p=10$. \label{fig:basisfunctions}}
\end{center}
\end{figure}

\myfig \ref{fig:basisfunctions} depicts two of the basis functions and their derivatives for $p = 10$ and $N = 80$. We can see that the magnitude of the basis functions may grow very large on the extended periodic domain, and the functions may also be highly oscillatory. Increasing $p$ will lead to greater oscillation while decreasing $p$ so that less points are used for the interpolation will lower the order of accuracy. Special care will need to be taken when computing integrals for the mass and stiffness matrices to resolve this behavior.
	
The proposed basis is not orthogonal, which means the matrices used in the DG framework will not be sparse, as they may be when using a standard nodal or modal basis. Therefore we will need to find an efficient way to calculate the inverse of the mass matrix and apply the mass and stiffness matrices in a PDE solver, especially in higher dimensions. Although it is not orthogonal, it is still a tensor-product nodal basis, which we will be able to take advantage of.	
	  
\section{Discontinuous Galerkin method with Fourier Continuation (FC-DG)}	
\label{sec:DGimplementation}
To illustrate the method, we first consider the discontinuous Galerkin formulation for the one dimensional transport equation. The techniques outlined are easily transferable to other  problems in one dimension. We then proceed to discuss extension to higher dimensions. 
		
\subsection{The Basic DG Method for the Transport Equation in One Dimension}
The scalar advection equation in a single space dimension can be written as
\begin{equation}
u_t + u_x = 0, \quad t>0, \quad x \in \Omega,  
\label{eqn:1dtransport}
\end{equation}
with initial condition $u(x,0) = f(x)$ and domain $\Omega$.
	
We discretize the domain $\Omega = [a,b]$ into $N_{\rm el}$ elements denoted $\Omega_k = [x_k,x_{k+1}]$, $k = 0,1,...,N_{\rm el} - 1$. Here $x_0 = a$ and $x_{N_{\rm el}} = b$. An approximation $u^h$ of the solution $u$ to \eqref{eqn:1dtransport} is then constructed element-wise as 
\begin{equation}
 u^h(x) = \sum_{j=0}^{N-1} \hat{u}^k_j \phi_j(x), \, x \in \Omega_k, \, k = 0,\ldots N_{\rm el}-1.
\label{eqn:solnapprox}
\end{equation}
Here $\phi_j$ are the basis functions from the test and trial space on each element, $N$ is the number of degrees of freedom on each element, and $\hat{u}_j$ are coefficients. Often a polynomial basis is used so $\phi_j \in P^{N-1}$, the space of polynomials of degree $N-1$, but to maintain generality we consider basis functions in some function space $V^h$. We of course intend to span this space by our Fourier continuation basis. 
	
To obtain the weak DG formulation, the approximation $u^h$ is substituted into \eqref{eqn:1dtransport}, multiplied by a test function $\phi_i$ from the same space $V^h$ as the basis functions and integrated over each element to obtain
\begin{equation}
0 = \int_{x_k}^{x_{k+1}} \phi_i u_t^h + \phi_i u_x^h dx.
\label{eqn:1dtransportelement}
\end{equation} 	

Applying integration by parts to the second term in the integrand of \eqref{eqn:1dtransportelement} and introducing a numerical flux term $u^*$ in the boundary terms results in the following element evolution equation  
\begin{equation}
	0 = \int_{x_k}^{x_{k+1}} \phi_i u_t^h dx - \int_{x_k}^{x_{k+1}} \frac{\partial \phi_i}{\partial x} u^h dx + \left[\phi_i  u^*\right]_{x_k}^{x_{k+1}}.
	\label{eqn:DG1dtransport}
\end{equation}
	
Different choices are possible for the numerical flux term $u^\ast = u^*(u^L,u^R)$, where $u^L$ and $u^R$ refer to the value of the approximation on the left and right side of a boundary respectively, in order to couple information between elements. 
In general, the numerical flux is required to be consistent, that is,  $u^*(u,u) = u$. 
For advective problems it is common to use an upwind flux. Here since the wave is traveling to the right the upwind flux becomes $u^* = u^L$. It can be shown that this choice of flux guarantees energy stability for the transport equation \cite{Hesthaven:2002ys}. 

Substituting the form of the approximation \eqref{eqn:solnapprox} into \eqref{eqn:DG1dtransport} gives 
\begin{equation}
\int_{x_k}^{x_{k+1}} \phi_i \sum_{j=0}^{N-1} \frac{\partial\hat{u}^k_j}{\partial t} \phi_j dx = \int_{x_k}^{x_{k+1}} \phi_i' \sum_{j=0}^{N-1} \hat{u}^k_j \phi_j dx - \left[\phi_i  u^*\right]_{x_k}^{x_{k+1}}. \label{eqn:DG1dtransport_weakform_originalelement}
\end{equation}	 
	
Here, to simplify the notation, we write derivatives with respect to $x$ using an apostrophe, i.e. $\partial \phi_i / \partial x = \phi_i'$. Requiring that (\ref{eqn:DG1dtransport_weakform_originalelement}) holds for each of the $N$ test functions $\phi_i$, $i = 0,...,N-1$ results in a system of $N$ equations on each element $\Omega_k$. The coupling between the element-wise systems is though the numerical flux. 

	
Since \eqref{eqn:DG1dtransport_weakform_originalelement} only differs between elements in the bounds of integration and the coefficients $\hat{u}_j$, but the same basis is used, we can write the weak formulation in a more generic way by mapping each element $\Omega_k = [x_k,x_{k+1}]$ to a reference element $[-1,1]$. We denote the spatial variable in the reference element by $z$. For the 1-D problem, this mapping is defined by the Jacobian $J_k = dx/dz = (x_{k+1}-x_k)/2$. Now we consider basis functions $\phi_j(z)$ on the reference element and the transformed weak formulation on element $\Omega_k$ can be written as
\begin{equation}
\int_{-1}^{1} \phi_i \sum_{j=0}^{N-1} \frac{\partial\hat{u}^k_j}{\partial t} \phi_j J_k dz = \int_{-1}^{1} \phi_i' \sum_{j=0}^{N-1} \hat{u}^k_j \phi_j dz - \left[\phi_i  u^*\right]_{-1}^{1}.
	\label{DG1dtransport_weakform}
\end{equation}	 
	
We define the \textit{mass matrix} $M$ and \textit{stiffness matrix} $S$ to have entries 
\begin{equation}
	M_{ij} = \int_{-1}^{1} \phi_i \phi_j J_k dz, \qquad S_{ij} = \int_{-1}^{1} \phi_i'\phi_j dz.
		\label{eqn:massstiffnessmtx}
\end{equation}
We also assemble the basis functions evaluated at the reference element boundaries into the \textit{lift matrices} $L_L$ and $L_R$ where $L$ and $R$ again denote the left and right boundary of the element.
\begin{equation}
		L_L = [\phi_0(-1),\phi_1(-1),...,\phi_{N-1}(-1)]^T,  L_R = [\phi_0(1),\phi_1(1),...,\phi_{N-1}(1)]^T.
\label{eqn:liftmtx}
\end{equation}
With this notation \eqref{DG1dtransport_weakform} can be written concisely in matrix-vector form as
\begin{equation*}
M\mathbf{\hat{u}}_t^k = S\mathbf{\hat{u}}^k + L_L \mathbf{u}^*_L - L_R \mathbf{u}^*_R.
\label{eqn:1dtransportsystem}
\end{equation*}
	
In a practical implementation, the flux terms $u^*$ will be computed first and then the time derivative $\mathbf{\hat{u}_t}$ can be found element by element by
\begin{equation}
	\mathbf{\hat{u}_t} = M^{-1} (S\mathbf{\hat{u}}+ L_L \mathbf{u^*_L} -  L_R \mathbf{u^*_R}).
	\label{eqn:u_t}
\end{equation}
Here $M^{-1} S$ is pre-computed and stored for efficiency.  
	
\subsection{Line-DG for Problems in Higher Dimensions}
\label{sec:LineDG}
	
To extend the DG formulation to higher dimensions, we use the Line-Based DG method described in \cite{PerssonLineDG2012}.  This scheme reduces connectivity of nodes within elements and thus increases the sparsity of the Jacobian matrix. To do this, 1-D DG solvers are used along each coordinate direction of the reference element. This circumvents the problem of inverting the dense higher dimensional mass matrix that would be generated using our new basis in the standard DG framework.

Where a standard nodal DG scheme would consider approximations in the space of 2-D polynomials, the line-based DG method considers each spatial derivative separately. As an example, consider solving the 2-D transport equation
\begin{equation*}
u_t + \alpha u_x + \beta u_y = 0, \qquad (x,y) \in \Omega,
\end{equation*}
on a rectangular element $\Omega$ with wavespeeds $\alpha, \beta > 0$. 
	
To obtain an approximation for $u_x$, a 1-D DG solver is applied along the $x$ direction at a number of fixed $y_j$, for $j = 0,...N-1$. After mapping to a reference element $[-1,1]$ as in the 1-D case, we define $u_j \in V^h([-1,1])$ to be the approximation function to $u$ that interpolates $u_{ij} = u(x_i,y_j)$, $i = 0,...N-1$ and $r_j$ to be the approximation to $u_x$ obtained from a 1-D DG formulation. The goal is to find $r_j \in V^h([-1,1])$ such that 
\begin{equation*}
	\int_{-1}^1 r_j(z) \cdot \phi(z) dz = \int_{-1}^1 \frac{du_j(z)}{dz} \cdot \phi(z) dz = -\int_{-1}^1 u_j(z) \cdot \frac{d\phi}{dz}dz + \left[u_j^* \cdot \phi\right]_{-1}^1.
\end{equation*}
where $u^*$ is given by some numerical flux function. Expanding $u_j$ and $r_j$ as an approximation in terms of the basis functions, we can substitute
\begin{align*}
	u_j(z) &= \sum_{i=0}^{N-1} \hat{u}_{ij}\phi_i(z), \\
	r_j(z) &= \sum_{i=0}^{N-1} \hat{r}_{ij} \phi_i(z),
\end{align*}
into the above formulation. The resulting system is equivalent to \eqref{eqn:u_t}. Solving for each $u_j$ gives the approximation to $u_x$ along the $x$-dimension at each grid point $(x_i,y_j)$, defined by $\hat{r}_{ij} = \hat{r}_{ij}^{(1)}$. The same procedure can be done to obtain the approximation for $u_y$ at each fixed $x_i$, $i = 0,\ldots,N-1$, for which we will denote the coefficients $\hat{r}_{ij}^{(2)}$. 
	
The final semi-discretized system for each $(u_{ij})_t$ is given by 
\begin{equation*}
\frac{d\hat{u}_{ij}}{dt} + \frac{1}{J}(\alpha \hat{r}_{ij}^{(1)} + \beta \hat{r}_{ij}^{(2)}) = 0.
\end{equation*} 
where $J$ is the determinant of the Jacobian mapping the physical element to the reference element $[-1,1]^2$. 
				
\subsection{Fourier Continuation Basis for DG}
For a nodal basis like the proposed FC basis, the coefficients $\hat{u}$ of the basis functions  are simply the values of the function evaluated at each node and $N-1$ is the number of grid points on the element. Hence, the approximation on the $k$-th element can be written as
\begin{equation}
u^h(x) = \sum_{l=0}^{N-1} u(x(z_l))\phi_l(z)J_k.
\label{eqn:nodalsolnapproximation}
\end{equation}	
	
The grid points in the reference element $z_l$ will be equidistant, ie. \mbox{$z_l = -1 + 2l/(N-1)$}. An equidistant grid is needed in order to use the FFT. 
The approximation of $u$ in 2-D can be written as
\begin{equation}
	u^h(x,y) = \sum_{i=0}^{N-1} \sum_{j=0}^{N-1} u_{ij} \phi_i(x) \phi_j(y)  = \sum_{i=0}^{N-1} \phi_i(x) \sum_{j=0}^{N-1} u_{ij} \phi_j(y),
\end{equation}
and similarly in higher dimension.
	
In order to evaluate the performance of our new basis, we will compare to a standard basis choice of Legendre polynomials $\phi_j = P_j(z)$ on Legendre-Gauss-Lobatto (LGL) nodes. Note that $N$ is still defined as the number of degrees of freedom, so we will use Legendre polynomials $P_j$ up to degree $q = N-1$. For this non-nodal basis, the initial data $u(x,0) = f(x)$ will first need to be expanded by an element-wise $L_2$-projection
	\begin{equation}
		M\hat{\textbf{u}} = 
		\begin{bmatrix}
		\int \phi_0 f(x)dx \\
		\int \phi_1 f(x)dx \\
		\vdots \\
		\int \phi_{N-1} f(x)dx 
		\end{bmatrix}.
		\label{eqn:L2approx}
	\end{equation}
	
Use of the Legendre polynomials results in a diagonal mass matrix due to the orthogonality of the basis. The choice of LGL nodes allows for integrals to be computed accurately up to degree $2N_{DG}-3$ by multiplying function values by pre-calculated LGL weights. In the new FC basis, we will not be able to use Gaussian quadrature since we have a uniform grid. We describe how we approximate the integrals in the next section.  

\subsection{Computing Integrals over the FC Basis}
One challenge is to compute the integrals exactly and efficiently using the new basis. At first this may appear difficult since we are using equidistant points and a non-polynomial basis. However, spectral interpolation can be used by zero-padding the FFT onto refined equidistant grids at a cost $\mathcal{O}(N \log N)$. We use the recent Gregory-type quadrature rules with interior weights 1 for equidistant grids introduced by Fornberg and Reeger \cite{Fornberg1Dquad2018} to obtain up to 16th order accuracy. It is also possible to exploit symmetry of the mass matrix to reduce the cost of its assembly. As we show in the experiments section, the condition number of the mass matrix is very small and does not depend on the number of gridpoints so the use of equidistant points does not affect the method adversely in terms of conditioning. Although not relevant for the linear problems considered here, it should be noted that for nonlinear problems where the integrals in the variational form of the flux must be evaluated at each timestep the oversampling will be expensive. 
	
Given the large magnitude of some of the basis functions on the extended domain as shown in \myfig \ref{fig:basisfunctions}, loss of accuracy may be experienced when taking the FFT for spectral interpolation or differentiating the basis functions. To deal with this, we generate the entries of the mass and stiffness matrices offline at high precision then convert them back to double precision for use in our PDE solvers. \texttt{MATLAB}'s multiprecision toolbox is used to do the high precision integral computations as it is compatible with \texttt{MATLAB}'s FFT.
	
\subsection{Time-stepping}
To step forward in time, a Taylor time stepping scheme is used. A Taylor series can be used to expand the solution around time $t+\delta t$ 
\begin{equation*}
u(t+\delta t) = u(t) + u_t \delta t + u_{tt}\frac{\delta t^2}{2!} + u_{ttt}\frac{\delta t^3}{3!} + ...
\end{equation*}

Given a semidiscretized system for the time derivative, $\mathbf{\hat{u}}_t = A\mathbf{\hat{u}}$, such as \eqref{eqn:u_t}, the discrete approximation to the time derivative terms on the right-hand side in the Taylor series above can be calculated sequentially as $\mathbf{\hat{u}}_t = A\mathbf{\hat{u}}$, $\mathbf{\hat{u}}_{tt} = A\mathbf{\hat{u}}_t$, etc. The number of terms taken in the Taylor series corresponds to the order of accuracy in the solution. 
	
For the centered and alternating flux the eigenvalues of the matrix $A$ will be purely imaginary and thus the timestepping method must have a stability domain that includes the imaginary axis. Taylor series methods with Taylor steps $N_t = 3,4,7,8,11,12,...$ have this property.  For our experiments we will use $N_t = 8$ unless otherwise specified.

\section{Numerical Experiments}\label{sec:numericalexperiments}
In this section we investigate the properties of the proposed basis and apply it to several test problems in 1-D and 2-D.
	
\subsection{Dispersive Properties of FC-DG}
As an initial test of the FC basis, we compute the numerical dispersion relation for the differentiation matrix resulting from \eqref{eqn:u_t}. For the linear transport equation in 1-D, the exact dispersion relation is given by $k = \alpha\omega$ where $k$ is the wavenumber, $\alpha$ is the wavespeed, which is positive, and $\omega$ is the frequency. A Bloch wave approach, as described in \cite{AinsworthDG}, is used to determine the numerical dispersion relation.
	
\begin{figure}[htbp]
\begin{center}
\includegraphics[width=0.50\textwidth]{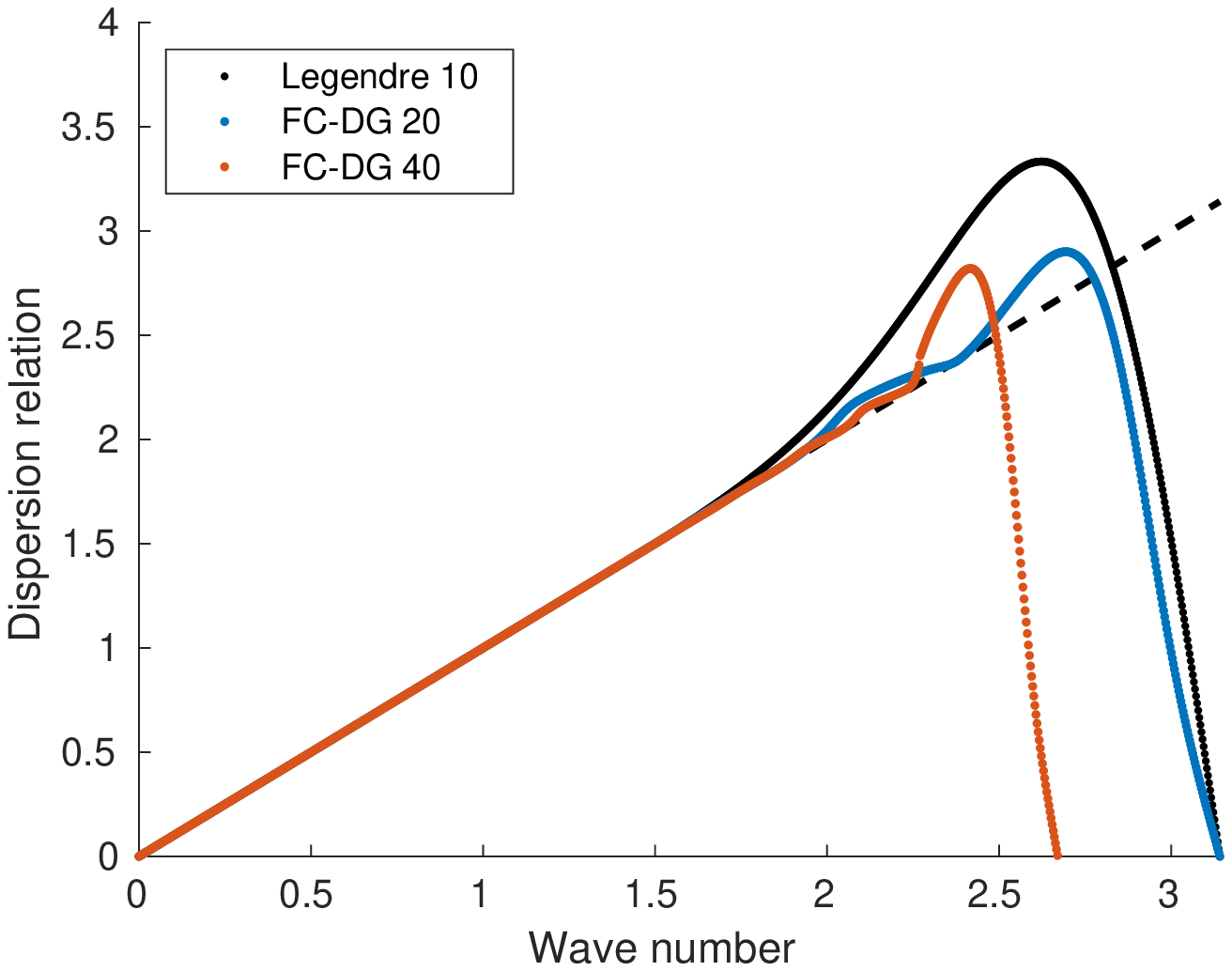}
\includegraphics[width=0.48\textwidth]{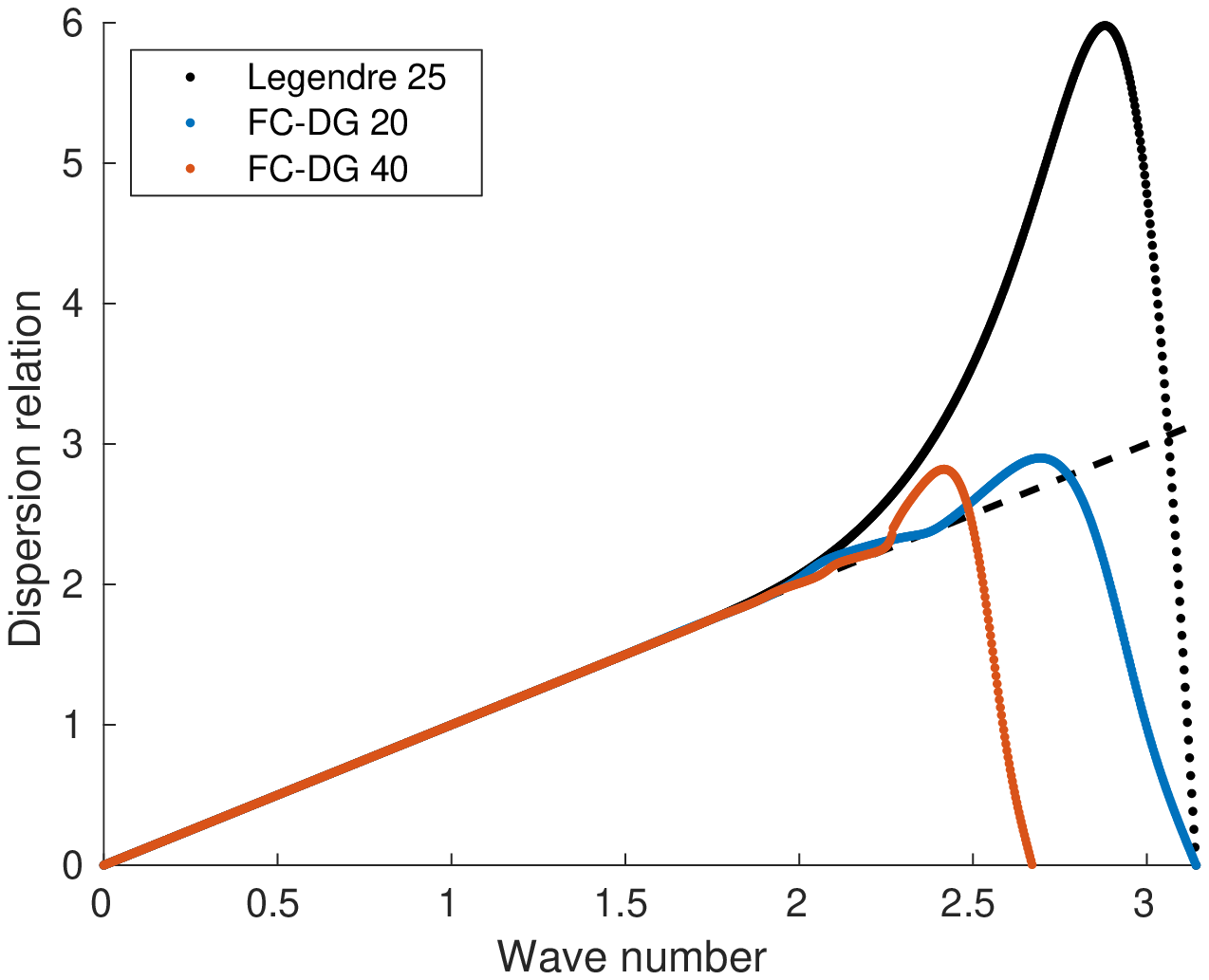}
\caption{Comparison of the dispersion relation for the 1-D transport problem for the Legendre basis and the FC basis. Left, 9th degree polynomials in the FC basis for $N = 20$ and $N = 40$ are compared to the Legendre basis of degree 10. Right, degree 25 Legendre polynomials have the same resolving power as the FC-basis.\label{fig:1dtransport_dispersion}}
		
\end{center}		
\end{figure}

The dispersion relation is obtained for the FC basis using degree 9 polynomials and $N = 20$, $40$ and $80$ gridpoints. This is compared to a 10th order and 20th order basis on Legendre-Gauss-Lobatto points. To compare, we look at the non-dimensional wave number $K = k/\Delta x$ and non-dimensional frequency $\Omega = \omega \Delta x/ \alpha$. \myfig \ref{fig:1dtransport_dispersion} depicts the normalized dispersion relations. Clearly the FC basis remains closer to the exact linear dispersion relation over larger wave numbers.
	
\begin{figure}[htbp]
\begin{center}
\includegraphics[width=0.7\textwidth]{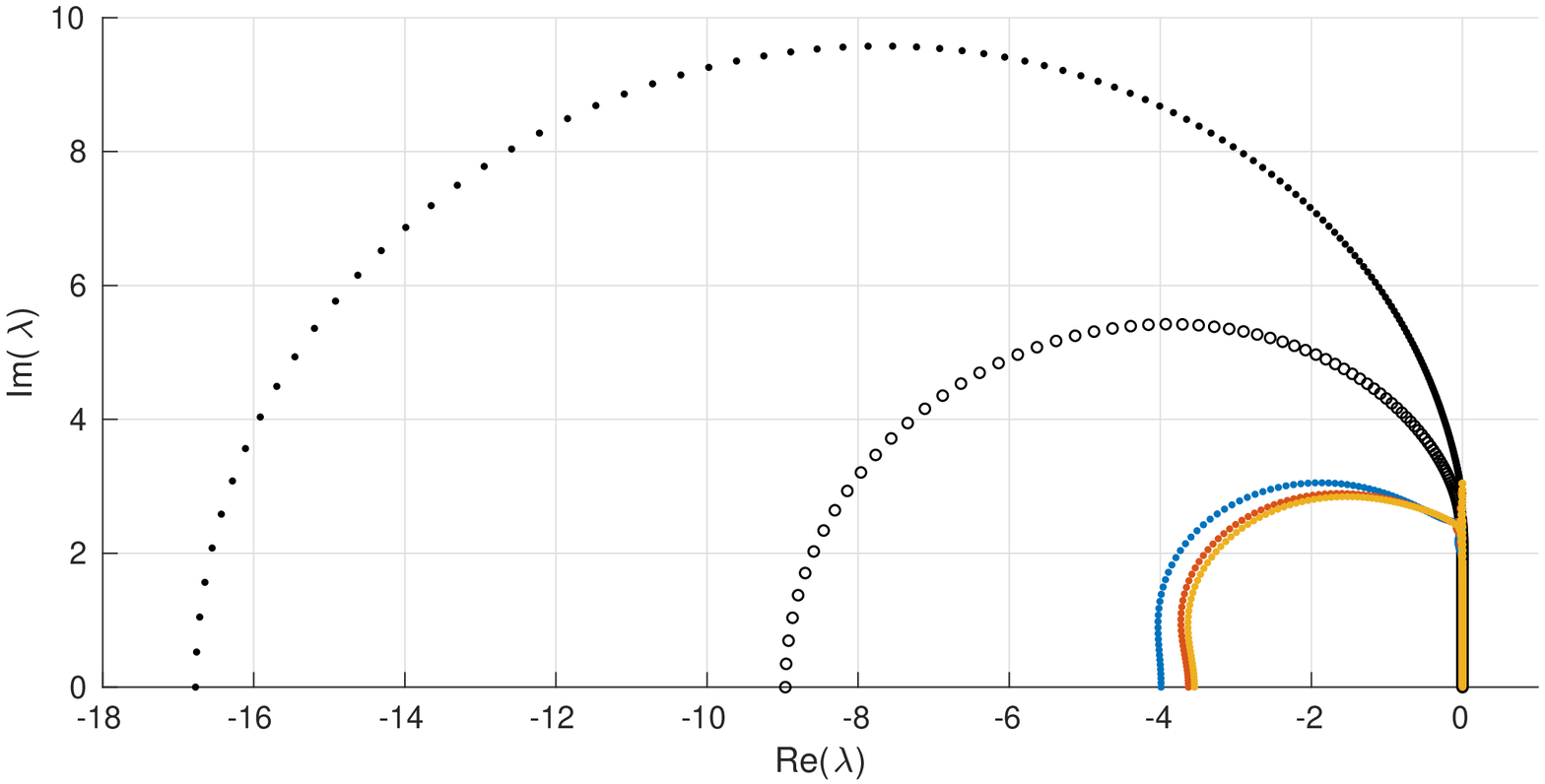}
\includegraphics[width=0.27\textwidth]{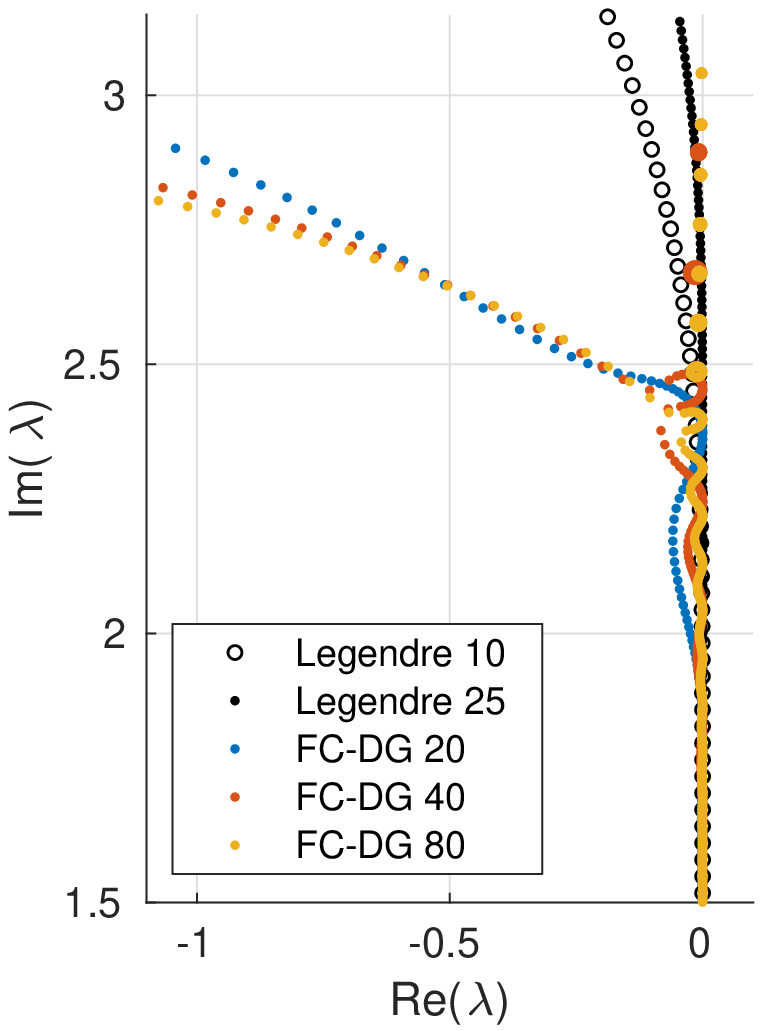}
\caption{Spectra of the differentiation matrix for the 1-D transport problem for the Legendre basis of degree 10 and 25, and the FC basis with $N = 20,40,80$. On the right, we zoom in to see the behavior at the imaginary axis.		\label{fig:1dtransport_spectrum}}
\end{center}
\end{figure}

\begin{figure}[htbp]
\begin{center}
\includegraphics[width=0.5\textwidth]{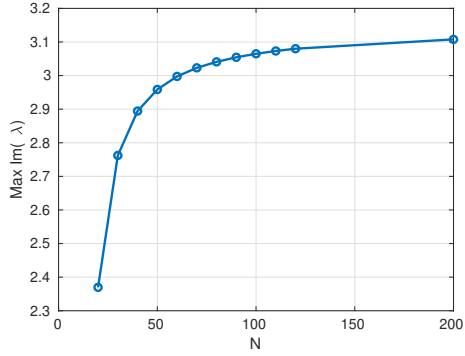}
\caption{Value of the largest eigenvalue along the imaginary axis as a function of degrees of freedom $N$. \label{fig:maxeigimag}}
\end{center}
\end{figure}

The full spectrum of the differentiation matrices are plotted in \myfig \ref{fig:1dtransport_spectrum}, scaled by the distance between nodes in an element and using 30 elements. Note that increasing the number of elements does not change the spectral radius, it just increases the number of eigenvalues lying along the curve. The spectral radius of the Legendre basis for $q = 20$ is more than 3 times larger than for the FC-basis, meaning that similarly larger timesteps can be taken when using the FC-basis than the Legendre basis at the same spatial resolution. Increasing $N$ does not significantly alter the spectral radius, which remains close to the rectangle $[-\pi,0] \times [\pi,\pi]$. \myfig \ref{fig:maxeigimag} depicts how the magnitude of the largest eigenvalue on the imaginary axis varies with $N$, further illustrating how this does not vary by much after $N$ is large enough and that it remains close to the limit of $\pi$. This illustrates the relationship of the FC basis to a pseudo-spectral discretization in terms of resolving power and time-stepping properties.	
	
As another metric, we compute the condition number of the mass matrix $M$ for the 1-D transport problem for various $N$ in Table \ref{table:conditionnumberM}. It can be seen that the condition number remains relatively constant for all $N$.  
	
\begin{table}[htbp]
\begin{center}
\begin{tabular}{lcccc}
\hline
$N$ & 20 & 40 & 80 & 200 \\
$\kappa(M)$ & 324.32 & 322.66 &  322.22  & 322.07\\
\hline
\end{tabular}	
\end{center}
\caption{Condition numbers for the mass matrix for various number of gridpoints on an element. As can be seen, the condition number is very robust with respect to changes in $N$. Here the FC basis using degree 9 polynomials was used.} 
\label{table:conditionnumberM}
\end{table}

\subsection{Long Time Errors}	
In this example we consider the transport equation in one dimension, \eqref{eqn:1dtransport} on a domain $x \in [-1,1]$ with periodic boundary conditions, and constant wave speed $\alpha = 1$. We use the upwind flux. Two cases for initial data cases are compared: $f_1(x) = \sin(10\pi x)$ and $f_2(x) = \exp(-50x^2)$, with respective analytic solutions $u_1(x,t) = \sin(10\pi(x-t))$ and $u_2(x,t) = \exp(-50(x-t)^2)$. We measure error in the $L^2$ norm as 
\begin{equation*}
	L^2\text{-error} = \left(\int_{-1}^{1}(u^h(x,t) - u(x,t))^2 dx\right)^{1/2}, \label{eqn:l2error}
\end{equation*}
where $u^h$ is the approximate solution and $u$ is the analytical solution. Given a discrete vector of values $u^h$ at nodes $x_j$, the integral \eqref{eqn:l2error} is computed using the standard trapezoidal rule. 
	
We compare the error using the new FC basis to a standard Legendre basis. For the FC basis, $N$ is the number of equidistant points used to construct the approximation in an element, and $N_{el}$ is the number of elements.  For the Legendre basis, $q$ is the degree of the polynomial approximation and $N_{el}$ is the number of elements. For the FC basis a CFL number of 0.2 is used, while 0.05 is used for the Legendre basis. 

\begin{figure}[]
\begin{center}	
\includegraphics[width=0.435\textwidth,trim={0.0cm 0.0cm 5.4cm 0.0cm},clip]{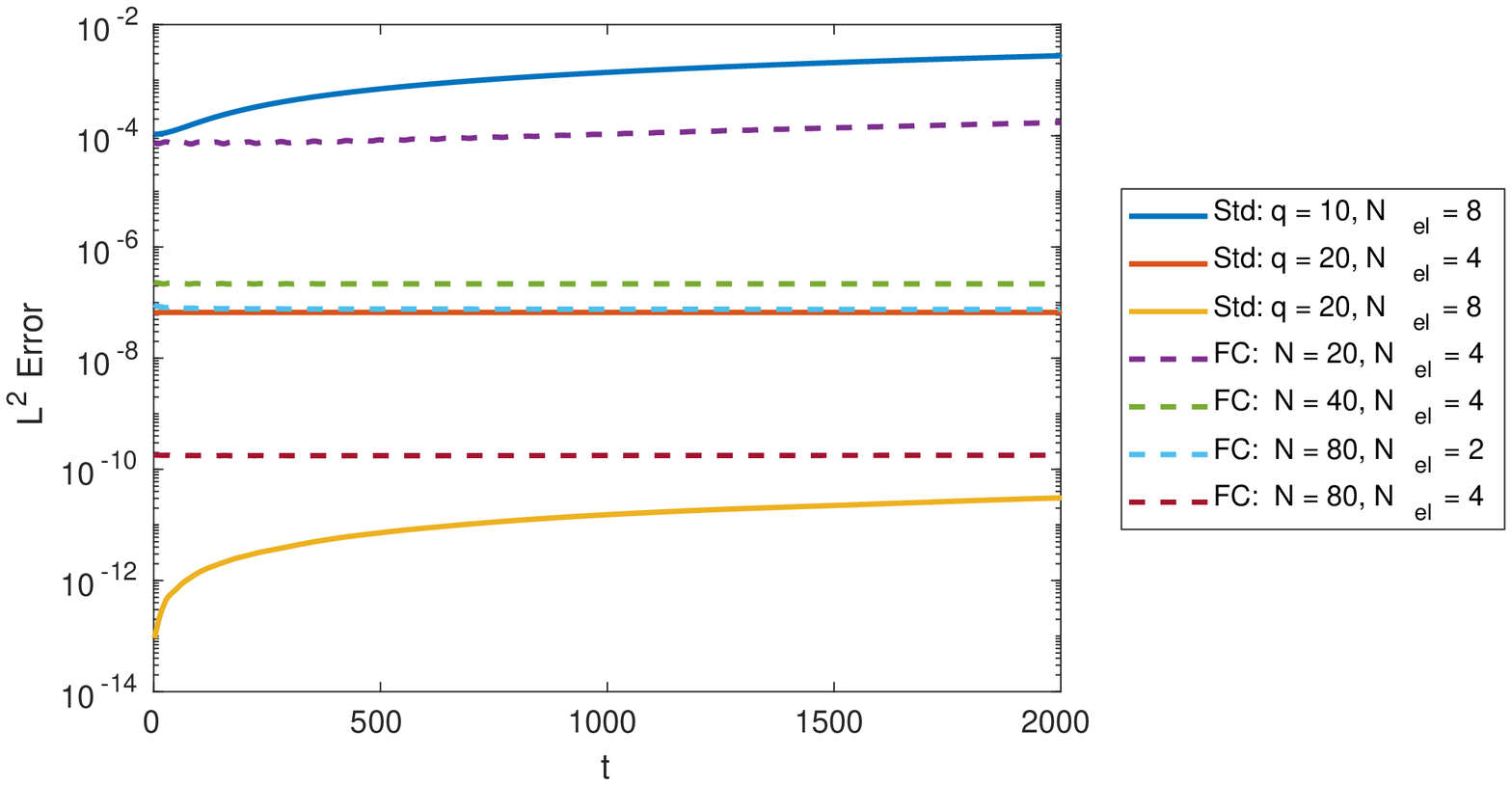}
\includegraphics[width=0.55\textwidth,trim={0.7cm 0.0cm 0.0cm 0.0cm},clip]{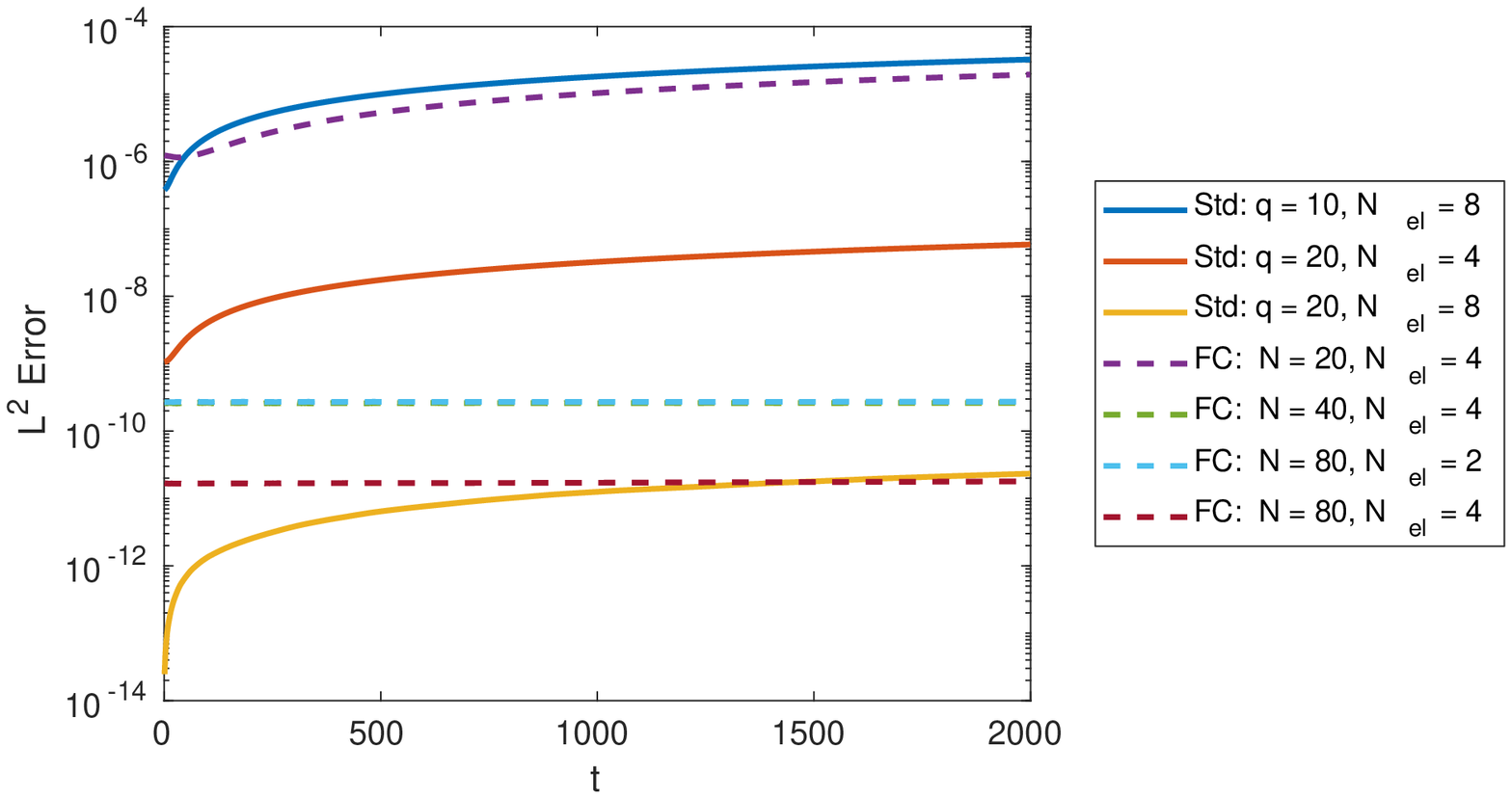}
\caption{Error over time for the 1-D transport of initial data $f(x)$. The results to the left are for $f(x) = \sin(10\pi x)$ and the results to the right are for $f(x) = e^{-50 x^2}$.		\label{fig:1dtransport_longrunerror}}
\end{center}
\end{figure}

\myfig \ref{fig:1dtransport_longrunerror} depicts the results for various degrees of freedom. It is evident that for large enough choices of $N$ and $N_{el}$, little to no dispersion is seen when using the FC basis. Although there is a trade-off in number of degrees of freedom, this becomes advantageous for problems that requires propagation of waves over many wavelengths in space or time.

\subsection{Investigation of Order of Accuracy}
\begin{figure}[htbp]
\begin{center}
\includegraphics[width=0.46\textwidth]{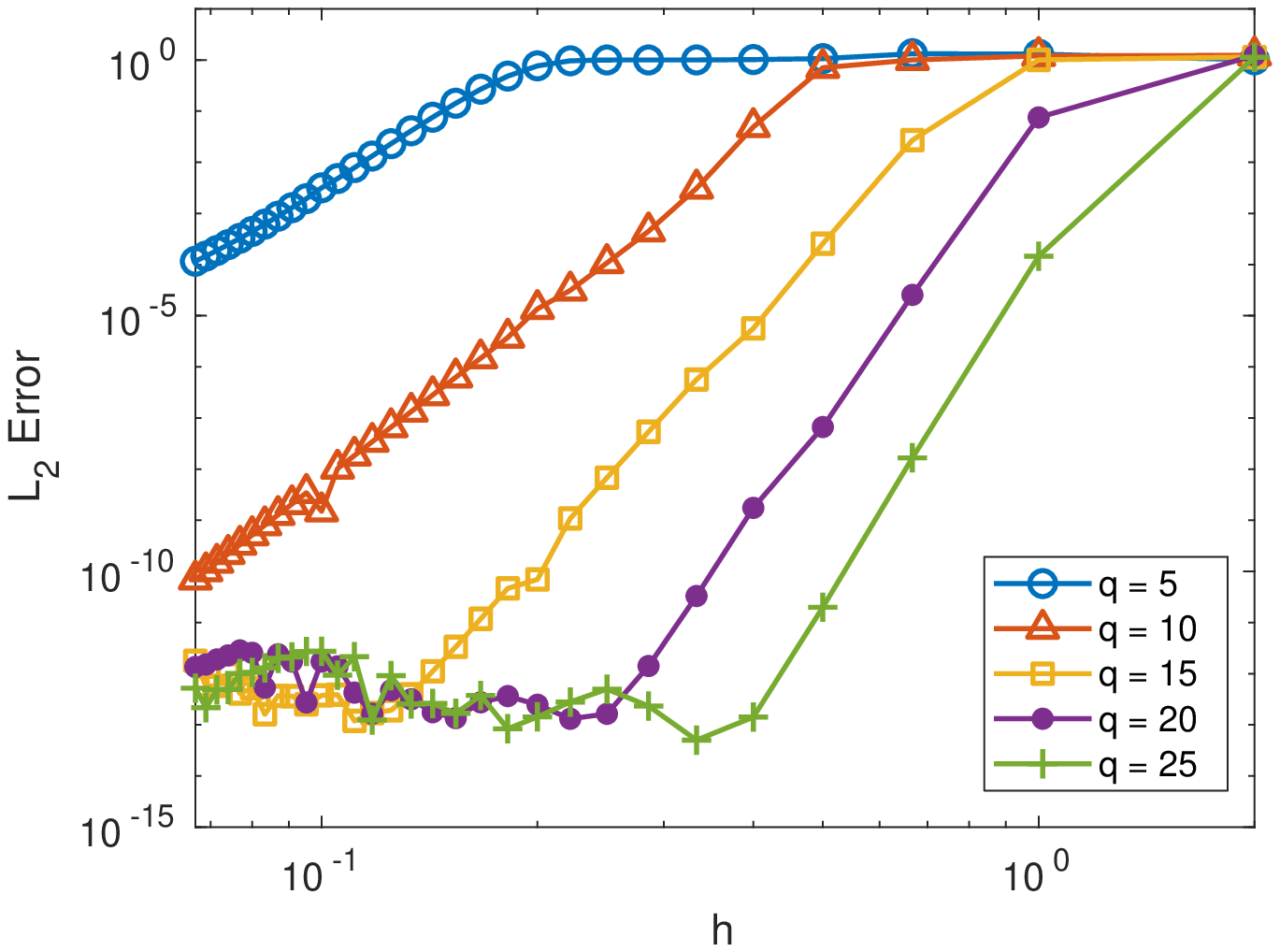}
\includegraphics[width=0.45\textwidth]{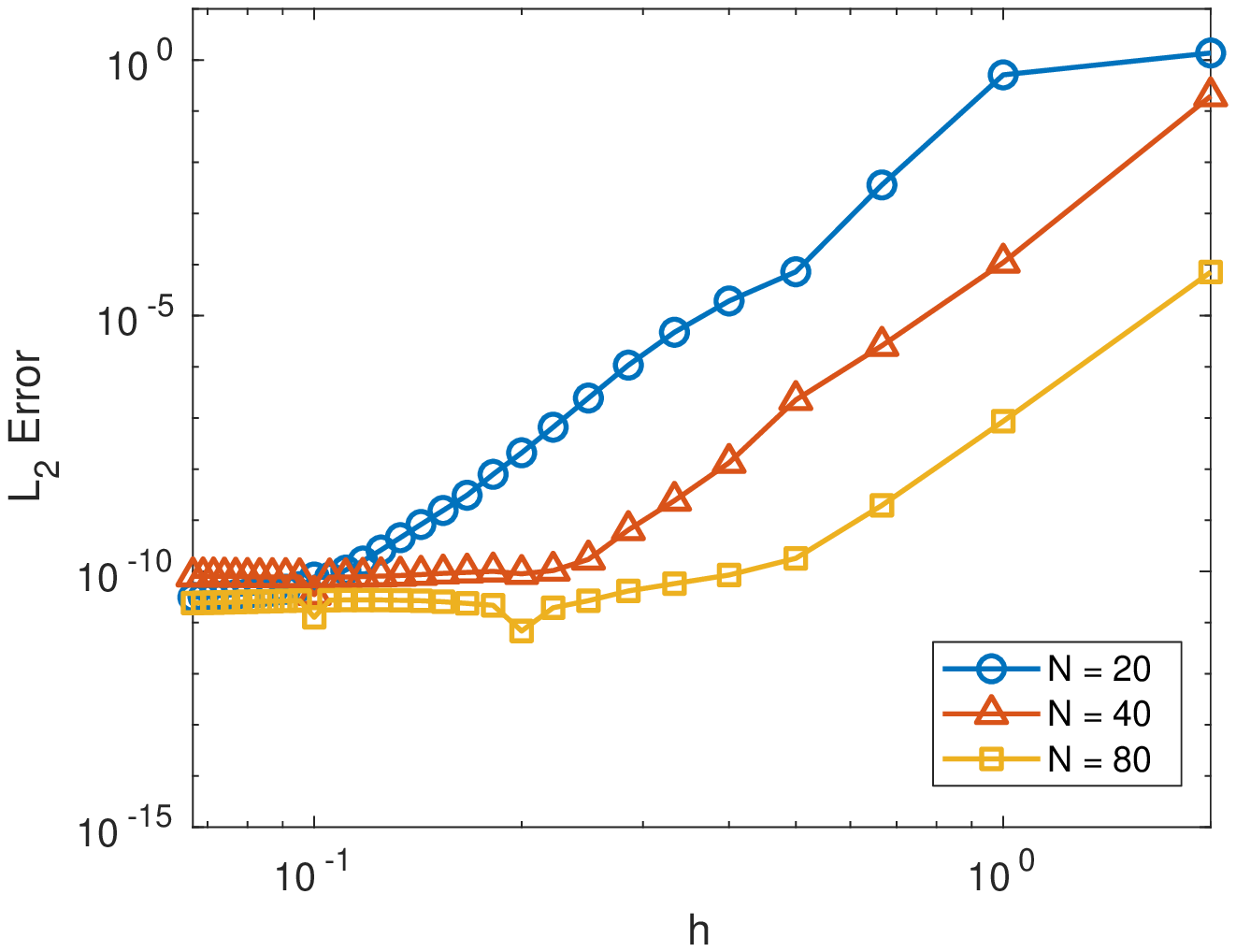}
\caption{Convergence of $L^2$ error for the 1-D transport equation with initial data $f(x) = \sin(10 \pi x)$. Left: Legendre basis. Right: FC basis.   \label{fig:1dtransport_conv}}
\end{center}
\end{figure}

Convergence of the approximation in the new basis is verified by measuring the error from the true solution to the 1-D transport equation at $T = 10$ for an increasing number of elements. Initial data is given by $f(x) = \sin(10\pi x)$, $x \in [-1,1]$, with wavespeed $\alpha = 1$ and upwind flux. Again, we use periodic boundary conditions. The convergence is plotted as a function of $h$, the length of each element, in \myfig \ref{fig:1dtransport_conv}. As expected, the convergence rates for the FC basis are approximately $p$, which is one order higher than the degree of interpolating polynomial used in the Fourier extension (computed convergence rates given in Table \ref{table:1dtransportconv}). Additionally, the rates of convergence do not depend on $N$. This can be compared to the standard Legendre basis using upwind fluxes, which has order of accuracy $q+1$ where $q$ is the degree of Legendre polynomial. 
	
Looking at \myfig \ref{fig:1dtransport_conv}, we can see that error in the approximation using the Legendre basis saturates due to machine precision around $10^{-12}$, but error using the FC basis saturates earlier around $10^{-10}$. 
	
\begin{figure}[htbp]
\begin{center}
\includegraphics[width=0.45\textwidth]{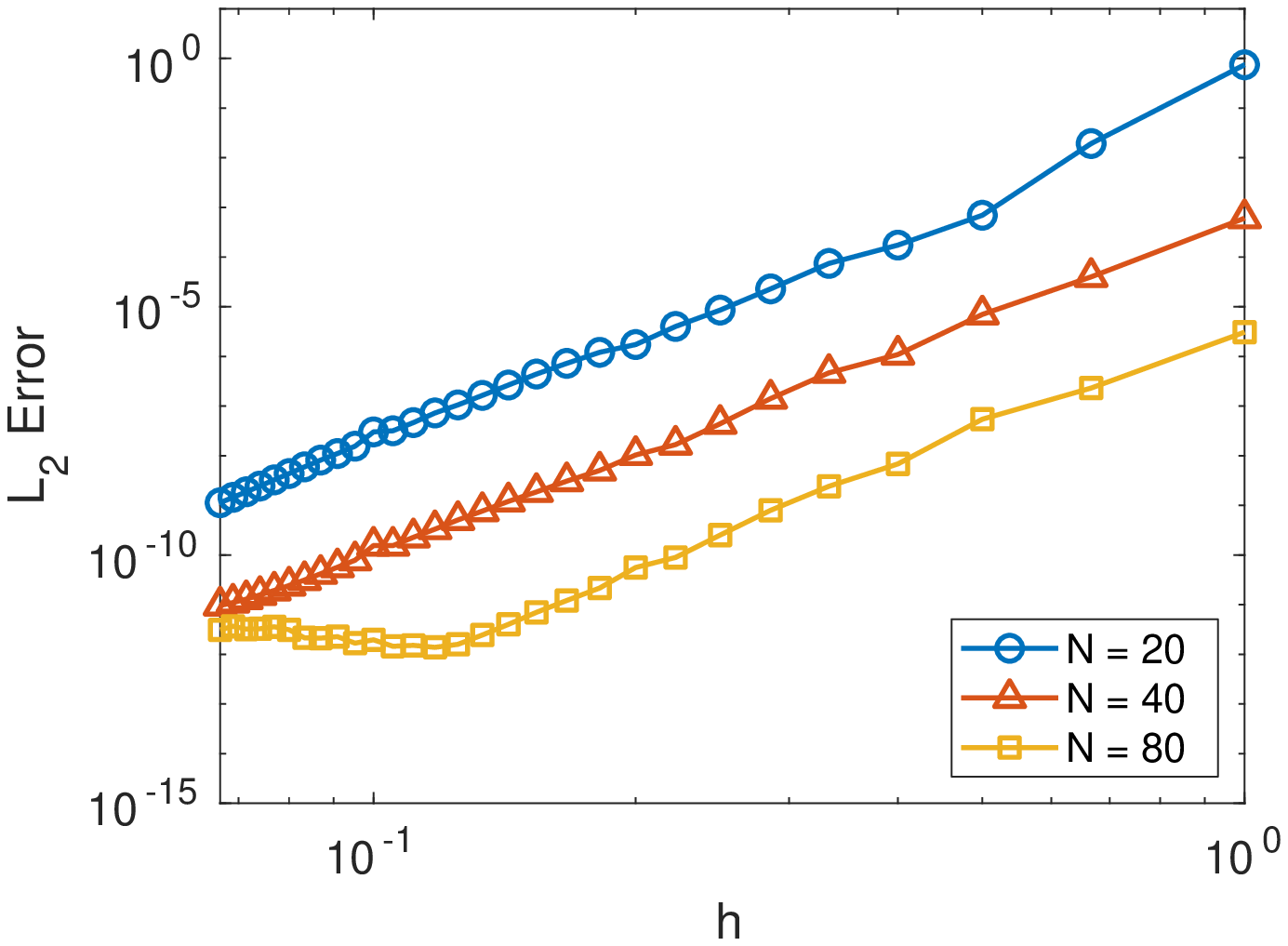}
\includegraphics[width=0.45\textwidth]{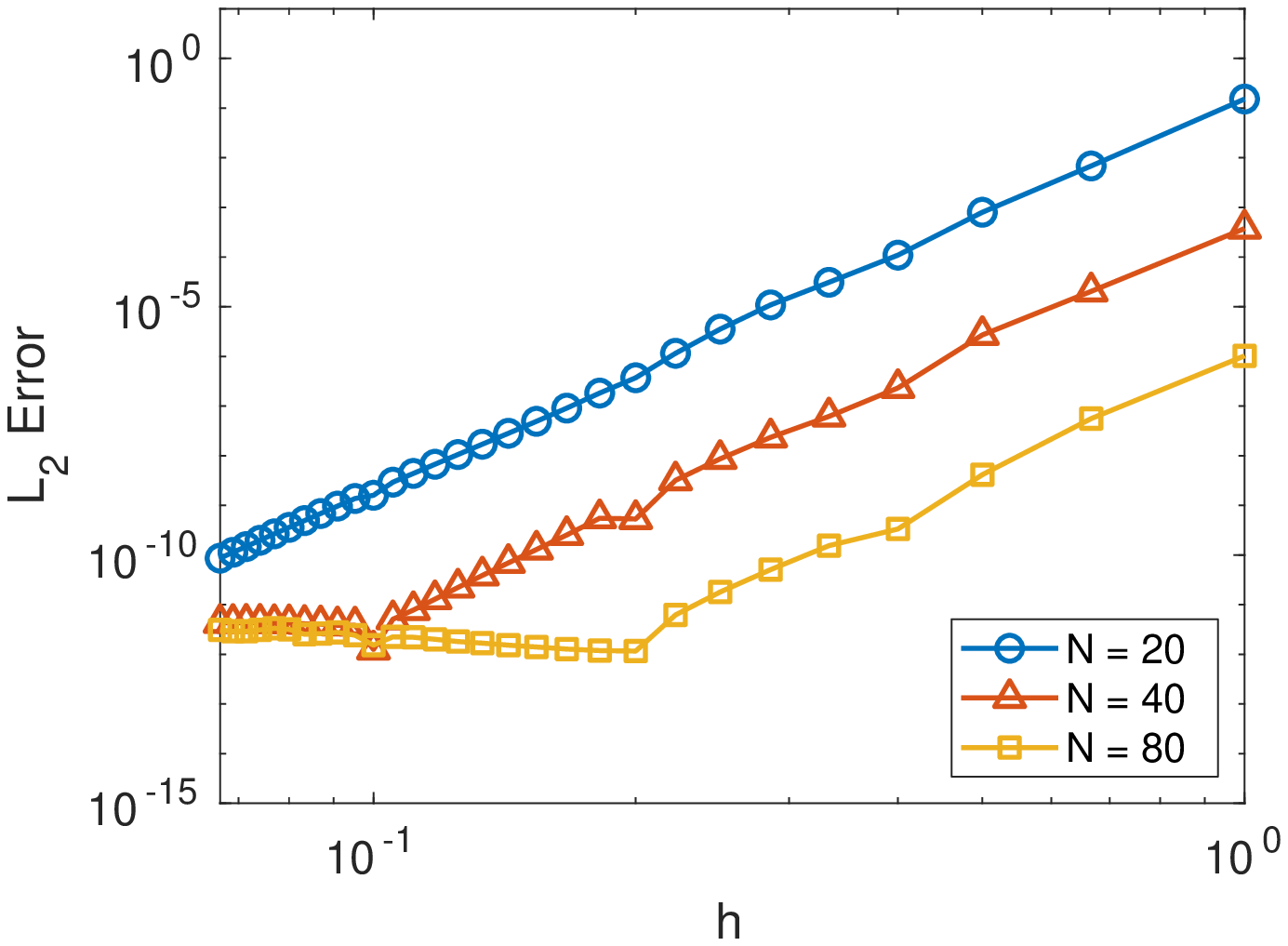}
\includegraphics[width=0.45\textwidth]{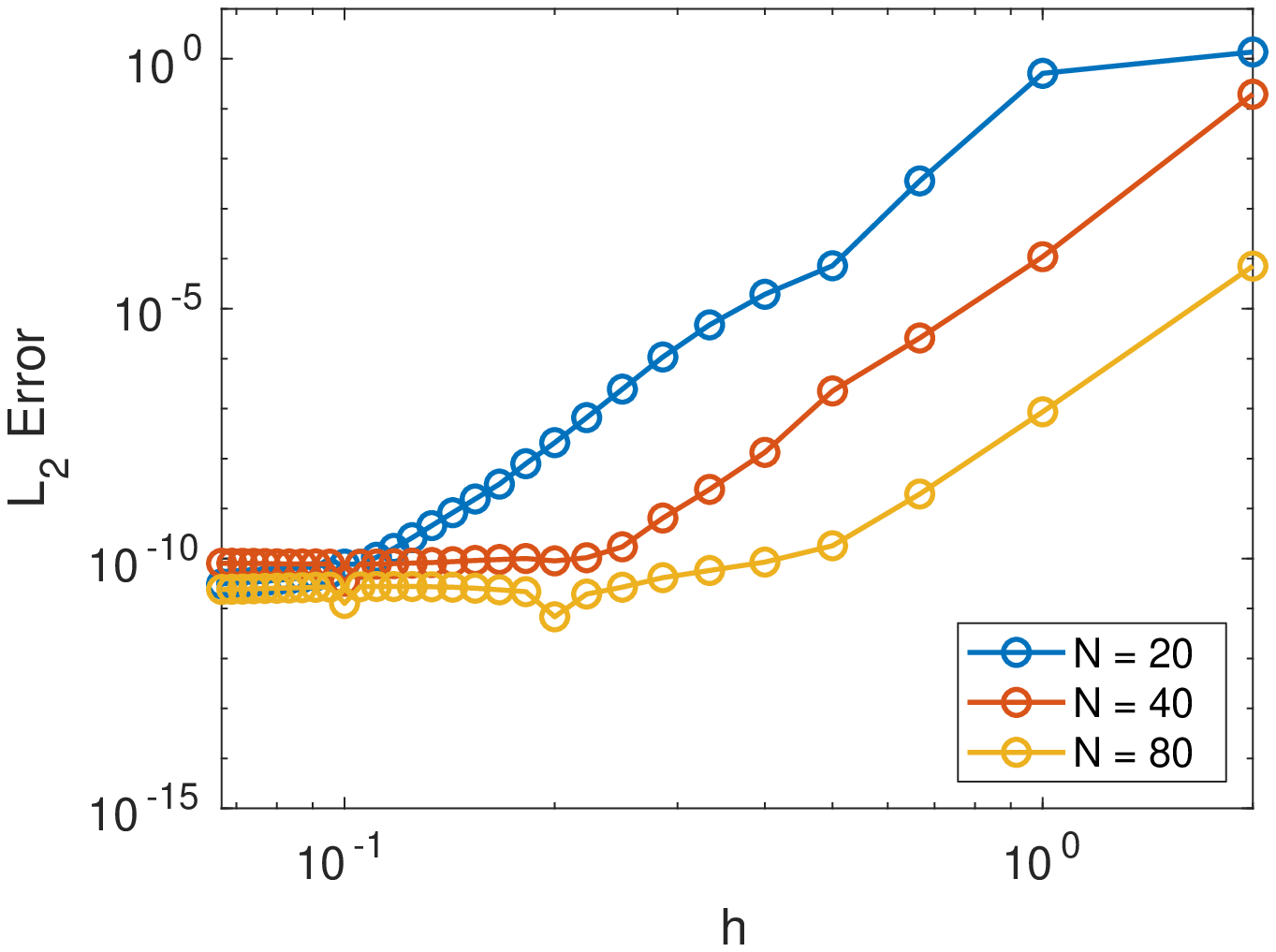}
\includegraphics[width=0.45\textwidth]{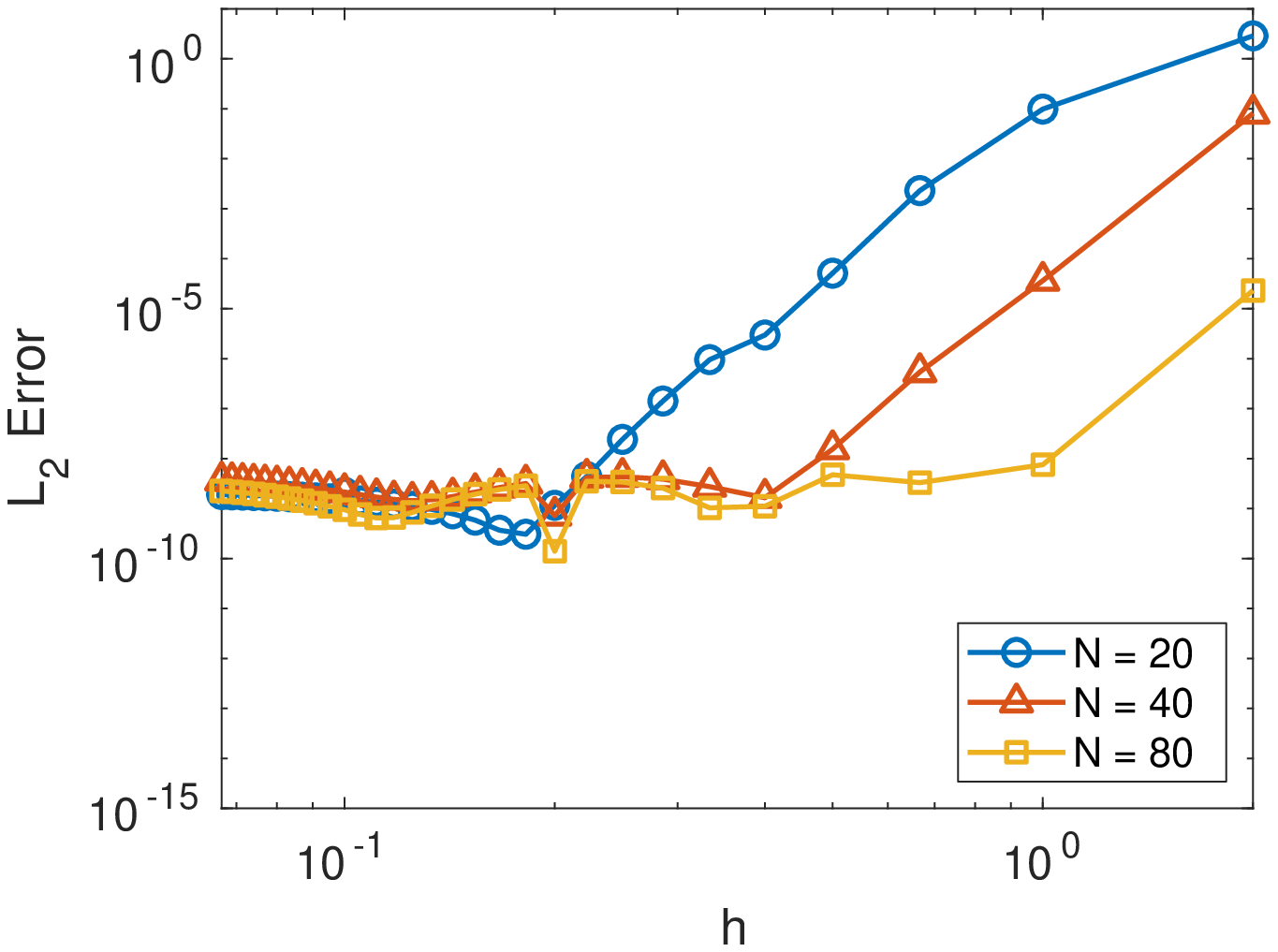}
\caption{Convergence of $L^2$ error for the 1-D transport equation using different degrees of polynomial in the Fourier continuation. From left to right, top to bottom, the degrees used are degree 6, 7, 9 and 11. 	\label{fig:1dtransport_deg_conv}}
\end{center}
\end{figure}

In \myfig \ref{fig:1dtransport_deg_conv}, we investigate the effect of the polynomial degree used in the Fourier continuation for the same 1-D transport problem described above. Note that the convergence for degree 9 polynomials is shown in both \myfig \ref{fig:1dtransport_conv} and \ref{fig:1dtransport_deg_conv}. Approximated convergence rates from a least squares fit are given in Table \ref{table:1dtransportconv}. Also noted is the approximate point at which the error saturates. The convergence rates remain approximately an order higher than the degree of polynomial until degree 11. Higher degree polynomials also appear to have a higher error saturation point.
	
\begin{table}[htbp]
\begin{center}
\begin{tabular}{lccc}
\hline
Number of basis functions  & 20 &  40 & 80 \\ 
\hline
Convergence rate deg 6 & 6.64 & 7.05 & 6.99 \\
Error saturation deg 6 & - & - & 1.59(-12) \\			
\hline
Convergence rate deg 7 & 7.50 & 8.10 & 8.01 \\
Error saturation deg 7 & - & 5.10(-12) & 1.17(-12) \\
\hline
Convergence rate deg 8 & 9.00 & 9.21 & 8.61 \\
Error saturation deg 8 & 2.00(-11) & 1.17(-11) & 2.55(-11) \\
\hline
Convergence rate deg 9 & 10.08 & 10.01 & 9.57 \\
Error saturation deg 9 & 1.06(-10) & 1.72(-10) & 1.80(-10)\\	
\hline		
Convergence rate deg 10 & 10.67 & 10.99 & 11.46 \\
Error saturation deg 10 & 3.85(-10) & 1.84(-10) & 9.81(-10)\\	
\hline	
Convergence rate deg 11 & 11.49 & 11.03 & 11.58\\
Error saturation deg 11 & 3.06(-10) & 1.65(-9) & 3.30(-9) \\	
\hline
\end{tabular}	
\end{center}
\caption{Convergence rates for the 1-D transport equation.\label{table:1dtransportconv}}
\end{table}

\subsection{Experiments with the Transport Equation in Two Dimensions}
Moving to higher dimensions, we solve the 2-D transport equation 
\begin{equation}
	u_t + \alpha u_x + \beta u_y = 0,
\end{equation}
on a structured grid of $N_{el}$ by $N_{el}$ elements in $[0,1]^2$ with initial data $u_0 = f(x,y)$ and periodic boundary conditions. Each element is discretized into $N$ by $N$ equidistant points that can be mapped to the reference element $[-1,1]^2$. The 2-D problem is discretized using Line-DG as described in Section \ref{sec:LineDG}. 
	
A convergence study is performed for increasing $N_{el}$ and various $N$ using upwind fluxes. The initial data is given by 
\[
f(x,y) = \sin(10\pi x) + \sin(10\pi y).
\] 
The error is measured after 1 cycle. 

\begin{figure}[htbp]
\begin{center}
\includegraphics[width=0.5\textwidth]{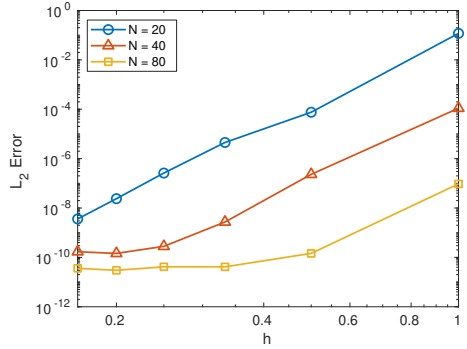}
\caption{Convergence of the $L^2$ error for 2-D transport.	\label{fig:2dtransportFCconv}}
\end{center}
\end{figure}

The results are shown in \myfig \ref{fig:2dtransportFCconv} and convergence rates are estimated in Table \ref{table:2dtransportconv}. Now we see convergence at a rate just higher than the degree of the interpolating polynomials in the extension. 

\begin{table}[htbp]
\begin{center}
\begin{tabular}{lccc}
	\hline
	Number of basis functions & 20 &  40 & 80 \\ 
	\hline
	Convergence rate & 9.46 & 9.43 & 9.35 \\
	\hline
\end{tabular}	
\caption{Convergence rates for the 2-D transport equation.	\label{table:2dtransportconv}}
\end{center}
\end{table}

\section{Application to Electromagnetic Waves}\label{sec:EMwaves}
Finally we consider applications to Maxwell's equations for describing electromagnetic waves, with particular consideration for behavior in optical materials. Maxwell's equations in a non-magnetic, non-conducting medium $\Omega \subset \mathbb{R}^d$, $d = 1,2,3$, $T> 0$, with no free charges, govern the dynamic evolution of the electric field $\mathbf{E}$ and the magnetic field $\mathbf{H}$, and can be written as 
\begin{subequations}\label{eq:max}
\begin{align}\label{eq:max1}
&\mu_0\ds\dd{t}{\mathbf{H}}+{\bf{\nabla}}\times \mathbf{E} = 0, \ \text{in} \  (0,T]\times \Omega,  \\[1.5ex]
\label{eq:max2}
&\epsilon_0\epsilon_\infty\ds\dd{t}{\mathbf{E}} +\epsilon_0\mathbf{J}-{\bf{\nabla}}\times \mathbf{H} = 0, \  \text{in} \ (0,T]\times \Omega, \\[1.5ex]
\label{eq:max3}
& {\bf{\nabla}}\cdot \mathbf{B}  = 0, \ {\bf{\nabla}}\cdot \mathbf{D} = 0, \ \text{in} \  (0,T]\times \Omega.
\end{align}
\end{subequations}
The electric flux density $\mathbf{D}$, and the magnetic induction $\mathbf{B}$, are related to the electric field and magnetic field, respectively, via the constitutive laws 
\begin{equation}
\label{eq:constD}
\mathbf{D} = \epsilon_0(\epsilon_\infty\mathbf{E}+\mathbf{P}), \ \ \mathbf{B} = \mu_0\mathbf{H},
\end{equation}
with the polarization current density, $\mathbf{J}$, defined as the time derivative of the macroscopic polarization, i.e. $\mathbf{J} = \dd{t}{\mathbf{P}}$. 
The parameter $\epsilon_{0}$ is the electric permittivity of free space, while $\mu_0$ is the magnetic permeability of free space. The term $\epsilon_\infty \mb{E}$ captures the linear instantaneous response of the material to the EM fields, with $\epsilon_{\infty}$ defined as the relative electric permittivity in the limit of infinite frequencies.

As an initial experiment, we consider dimension $d = 2$ and take $\mathbf{J} = \mathbf{P} = 0$. We can write the simplified evolution equations component-wise as 
\begin{subequations}
\begin{align}
\mu_0 \frac{\partial{H^z}}{\partial t} &= -\frac{\partial{E^y}}{\partial x}+\frac{\partial{E^x}}{\partial y}, \\
\epsilon_0 \epsilon_\infty \frac{\partial{E^x}}{\partial t} &= \frac{\partial{H^z}}{\partial y },\\
\epsilon_0 \epsilon_\infty \frac{\partial{E^y}}{\partial t} &= -\frac{\partial{H^z}}{\partial x},
\end{align}
\label{eqn:maxwell2d}
\end{subequations}
where $E^x = E^x(x,y,t)$ and $E^y = E^y(x,y,t)$ are the $x$ and $y$ components of the electric field , $H^z = H^z(x,y,t)$ is the magnetic field in the $z$-direction. $E^x$ and $E^y$ are constrained to be $0$ at tangential boundaries. The initial magnetic field is prescribed to be a function $H_z(x,y,0) = f(x,y)$ and the initial $E^x$ and $E^y$ fields are both set to be zero. 
The tangential boundary condition for $E^x$ and $E^y$ is implemented by setting $E^+ = -E^-$ at those exterior boundaries. All other exterior boundaries are set so $E^+ = E^-$, resulting in a first-order Neumann boundary coundition.

Using normalized parameters $\mu_0 = \epsilon_0 = \epsilon_\infty = 1$, initial condition $f(x,y) = \sin(5x)\sin(5y)$, and a domain of $\Omega = [-3\pi/2,3\pi/2]^2$, the numerical model is evolved for one period and the resulting $H^z$ field compared to the analytical solution 
\[
H^z(x,y,t) = \sin(5x)\sin(5y)\cos(5\sqrt{2}t).
\] 

\begin{figure}[htbp]
	\centering
	\includegraphics[width=0.5\textwidth]{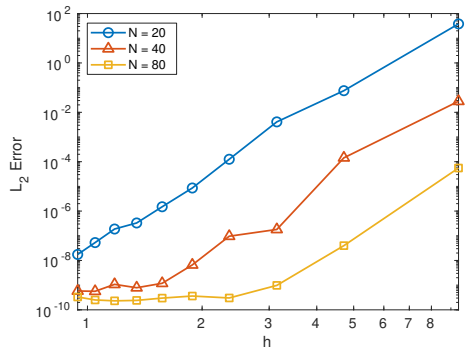}
	\caption{Convergence of the $L^2$ error for 2-D Maxwell's equations.}
	\label{fig:2dmaxwellconv}
\end{figure}

\myfig \ref{fig:2dmaxwellconv} shows the resulting error under refinement for various $N$, grid points per element, and Table \ref{table:2dmaxwellconv} gives estimated convergence rates using both centered fluxes and alternating fluxes. Similar to the results for the 2-D transport equation, the order of convergence is seen to be slightly higher to a degree higher than the order of polynomials used in the Fourier extension ($p = 9$).

\begin{table}[htbp]
	\centering
	\begin{tabular}{lccc}
		\hline
		Number of basis functions  & 20 &  40 & 80 \\ 
		\hline
		Convergence rate C.-Flux & 9.49 & 9.68 & 10.02 \\
		\hline		
		Convergence rate A.-Flux & 9.42 & 9.43 & 10.18 \\
		\hline	
	\end{tabular}	
	\caption{Convergence rates for the 2-D Maxwell's equations for centered fluxes (C.-Flux) and alternating fluxes (A.-Flux).}
	\label{table:2dmaxwellconv}
\end{table}

To demonstrate the capability of handling more complex solutions, we again consider Maxwell's equations \eqref{eqn:maxwell2d}, but introduce a forcing term $f(x,y,t)$ as
\begin{subequations}
\begin{align}
\mu_0 \frac{\partial{H^z}}{\partial t} &= -\frac{\partial{E^y}}{\partial x}+\frac{\partial{E^x}}{\partial y}, \\
\epsilon_0 \epsilon_\infty \frac{\partial{E^x}}{\partial t} &= \frac{\partial{H^z}}{\partial y } + f(x,y,t)(y-y_0),\\
\epsilon_0 \epsilon_\infty \frac{\partial{E^y}}{\partial t} &= -\frac{\partial{H^z}}{\partial x } + f(x,y,t)(x-x_0),
\end{align}
\label{eqn:maxwellforcing}
\end{subequations}	
where $(x_0,y_0)$ is a given source point in the domain. To implement this within our numerical scheme, the forcing term is added in point-wise at each timestep.

\begin{figure}[htbp]
	\includegraphics[width=\textwidth]{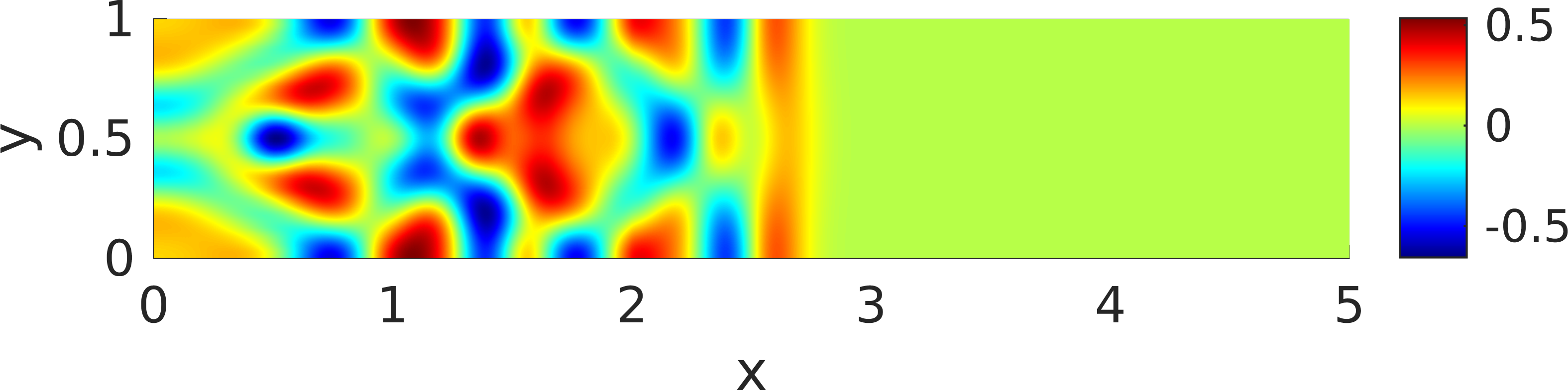}
	\includegraphics[width=\textwidth]{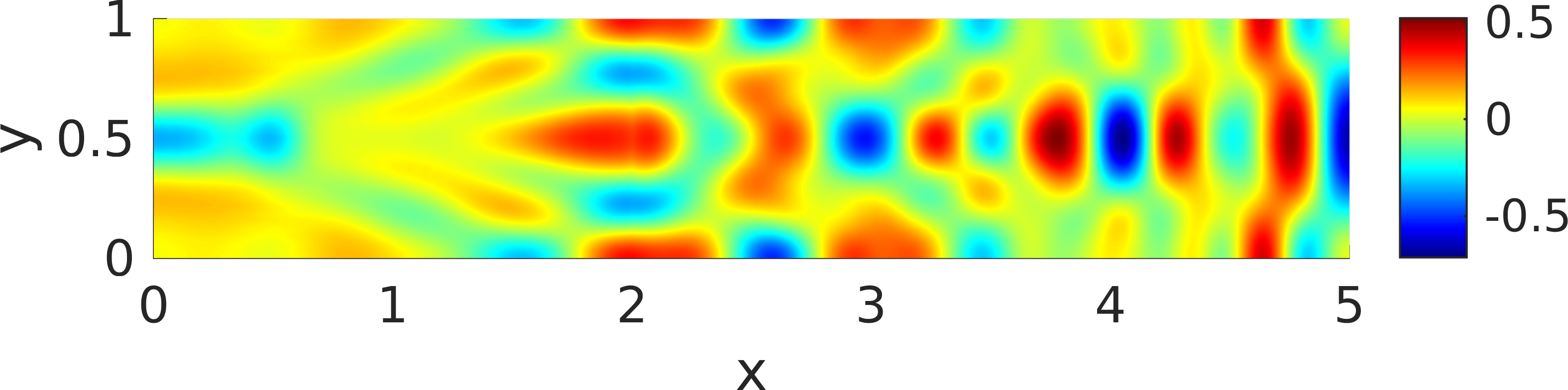}
	\includegraphics[width=\textwidth]{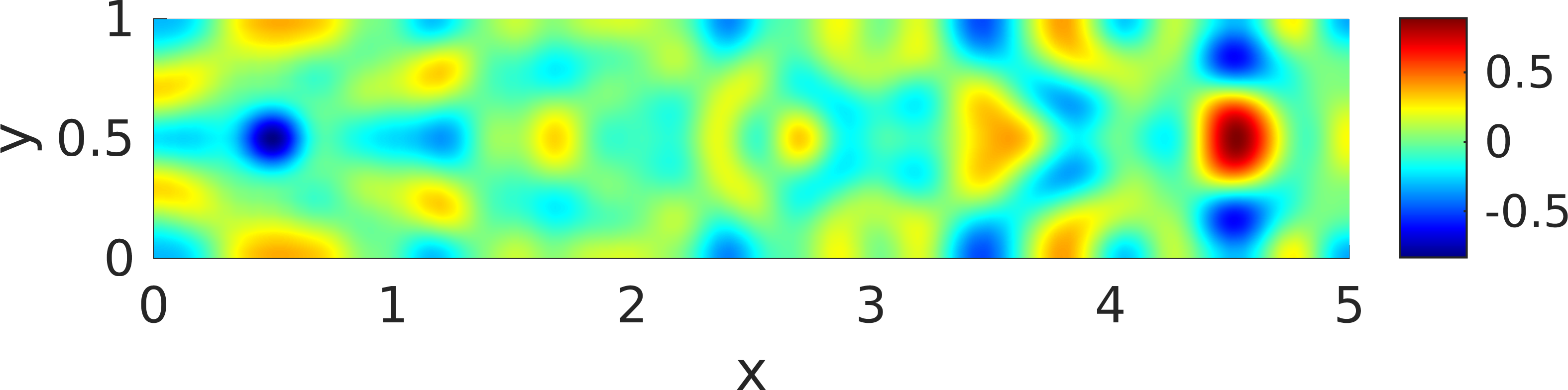}
	\caption{Snapshots of the magnetic field subject to forcing function.		\label{fig:maxwellsnapshots}}
\end{figure}

In this experiment, initial data is given by $H^z(x,y,0) = 0$ on $[0,1]\times[0,5]$ with parameters $\mu_0 = \epsilon_0 = \epsilon_\infty = 1$. The forcing term is given by 
\[
f(x,y,t) = 50^2 \sin(100 t)\exp(-36((x-x_0)^2+(y-y_0)^2)),
\]
 with $(x_0,y_0) = (0.5,0.5)$.  The number of elements used is 10 in the $x$ direction and 2 in the $y$ direction, with $N = 40$ on each element. The solution for the $H^z$ field is depicted at $T = 2, 5, 50$ in \myfig \ref{fig:maxwellsnapshots}.

In order to add in non-zero polarization terms $\mb{J}$ and $\mb{P}$, we use the ADE approach as in \cite{IEEE_Duffing} and append a system of ODEs describing the nonlinear relationship between the macroscopic polarization vector field ${\mb P}$ and the electric field ${\mb E}$ to Maxwell's equations.
The macroscopic  {\em (electric) polarization} $\mathbf{P}$ includes both linear and nonlinear effects, and is related to the electric field through different mechanisms depending on the optical phenomenon under consideration.  In this work we consider what is known as {\em general Maxwell-Duffing dispersive models.} 
\begin{figure}[htb]
\begin{center}
\includegraphics[width=0.3\textwidth]{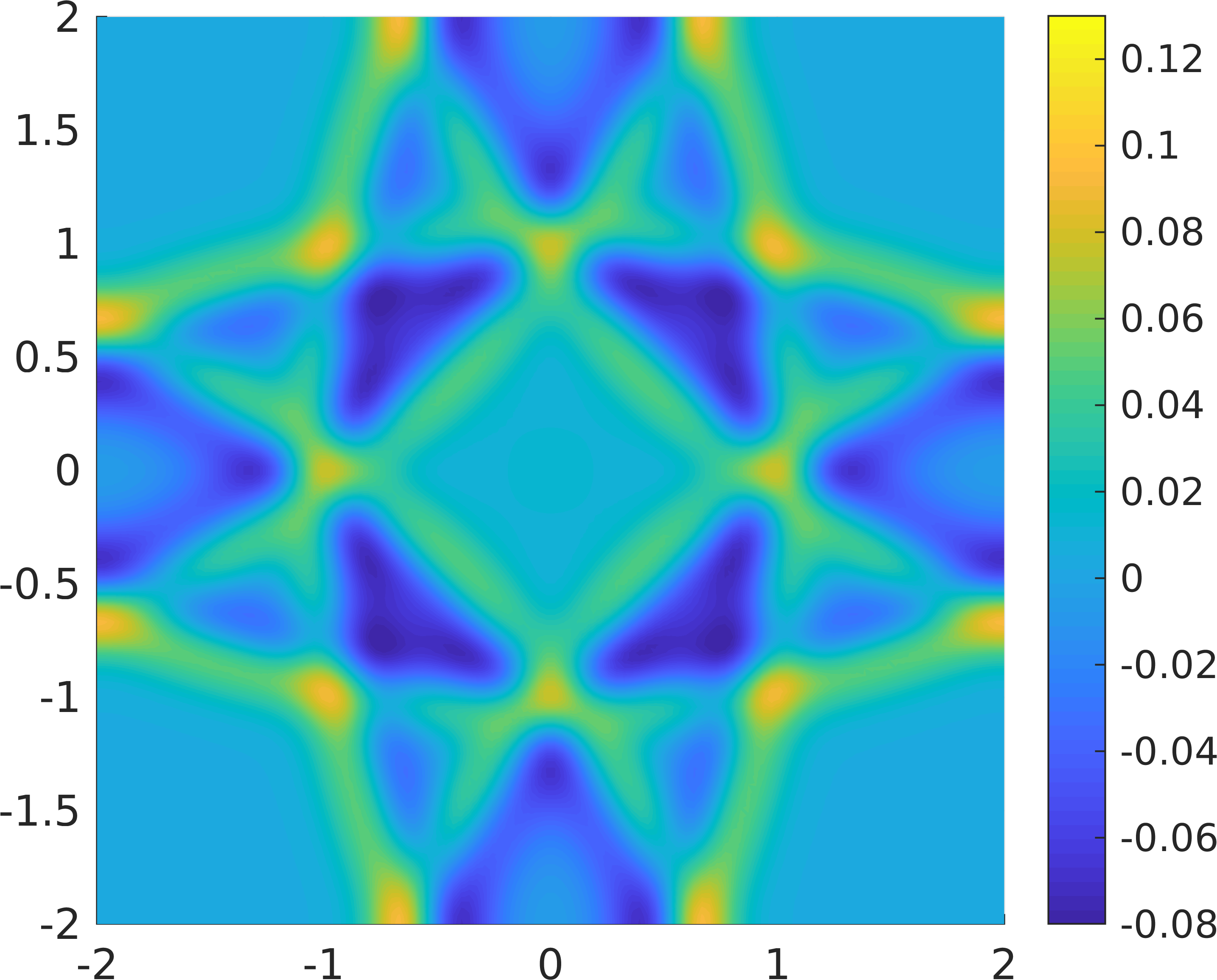}
\includegraphics[width=0.3\textwidth]{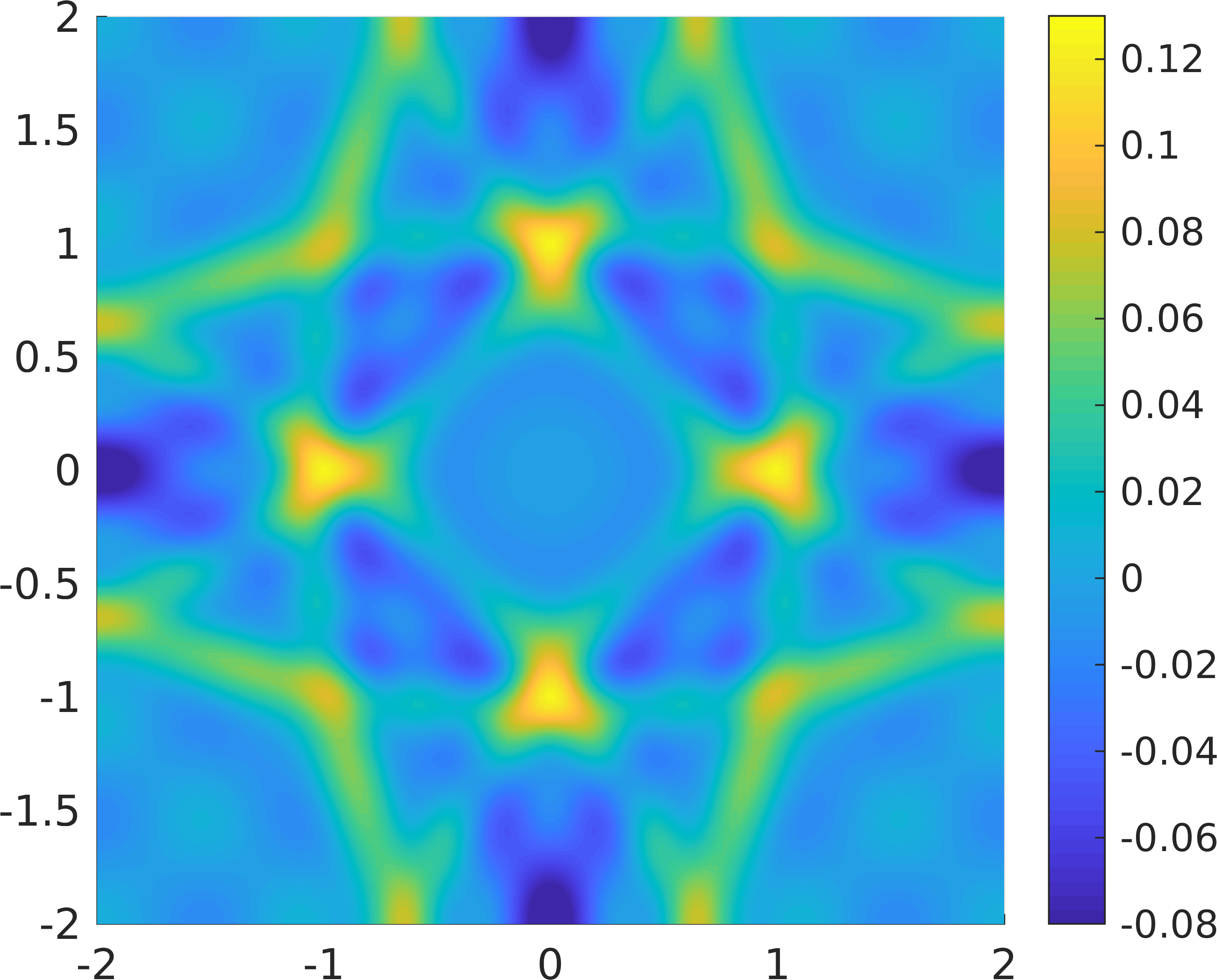}
\includegraphics[width=0.3\textwidth]{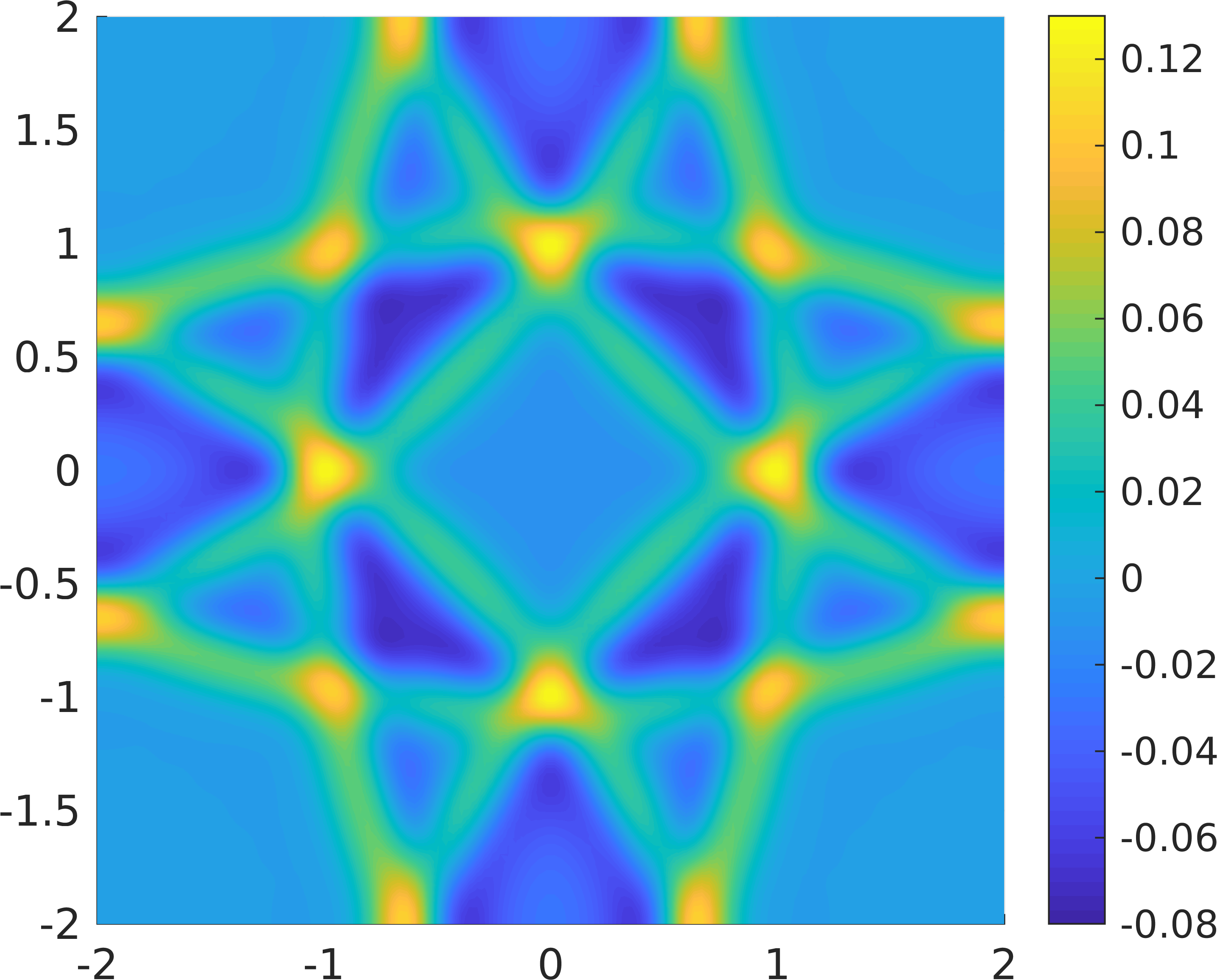}
\caption{Three snapshots of the solution to the Duffing model with all material parameters set to unity except $\omega_0$ which is set to 1, 100 and 1000 from left to right. Displayed is the $H_z$ field. \label{fig:duff}}
\end{center}		
\end{figure}

The Duffing equation for the electric polarization, models high order effects by including both nonlinearity and dispersion, and can be written in a general form as 
 \begin{equation}
  \label{eq:duff}
\frac{\partial^2 \mathbf{P}}{\partial t^2}+\frac{1}{\tau}\frac{\partial\mathbf{P}}{\partial t} +\omega_0^2\mathbf{P}F(\mathbf{P})=\omega_p^2\mb{E},
\end{equation}
with a range of possible choices for $F(\mathbf{P})$ \cite{}. Here $\omega_0$ and $\omega_p$ are the resonance and plasma frequencies of the medium, respectively, and $\tau^{-1}$ is a damping constant.  We will consider an Nth order polynomial model for the Duffing equation, given as 
\begin{equation}
  F(\mathbf{P})=F_{\rm PMD}(\mathbf{P}) : =  \sum_{l=0}^{N_{\rm PMD}} \lambda_{2l} | \mathbf{P} |^{2l},
  \label{eq:PMD}
\end{equation}
with $N_{\rm PMD} \in \mathbb{N}, N_{\rm PMD} \geq 1$. We refer to the system of equations obtained by adding \eqref{eq:duff} and \eqref{eq:PMD} to \eqref{eq:max} as the  {\em Nth Order Polynomial Maxwell-Duffing (PMD) model}.

We note that if $F(\mathbf P) = 1$, the Duffing model reduces to the linear Lorentz dispersive model. A sample computation using this model in 2-D with $N_{PMD} = 1$ and all material parameters but $\omega_0$ set to unity can be found in \myfig \ref{fig:duff}.

\section{Conclusion}\label{sec:conclusion}
This paper has presented a new method, the Fourier continuation - discontinuous Galerkin method, constructed by utilizing the discrete Fourier extension from \cite{AlbinDPE2014} as a basis in the traditional discontinuous Galerkin framework \cite{cockburn1989tvb,Hesthaven:2002ys}. We demonstrated through numerical experiments that our method has good dispersive error properties and that these properties translate to accurate propagation of waves over many wavelengths. Our method also admits larger timesteps than traditional polynomial based DG methods. 

The main drawbacks of the method are: 1.) the reliance on oversampling (through spectral FFT interpolation) which makes the assembly process more expensive than for methods that collocate the degrees of freedom and the quadrature nodes, and 2.) the non-orthogonality of the basis resulting in dense mass matrices.  

We believe that both of these drawbacks may be overcome or at least mitigated in future research. Here we have exclusively focused on hyperbolic problems, but note that the framework would also be possible to extend to elliptic problems. There block preconditioning with traditional FC solvers could prove fruitful. 

\section{Conflict of Interest Statement}
On behalf of all authors, the corresponding author states that there is no conflict of interest.

\bibliographystyle{plain}
\bibliography{references,appelo}
\end{document}